\documentclass[10pt, twoside]{scrartcl}
\usepackage{etex}
\usepackage[a4paper, top=88pt, bottom=88pt, left=72pt, right=72pt, headsep=16pt, footskip=28pt]{geometry}
\usepackage{amsmath, amssymb, amsthm, xcolor, tabularx}
\usepackage[hyphens]{url}
\usepackage[colorlinks, citecolor=, linkcolor=, urlcolor=]{hyperref}
\usepackage[nameinlink]{cleveref}
\usepackage{enumitem}
\usepackage{mleftright} \mleftright
\renewcommand{\qedsymbol}{$\blacksquare$}

\usepackage{scrlayer-scrpage, xhfill}
\automark[section]{section}
\setkomafont{pageheadfoot}{\normalcolor\sffamily}
\setkomafont{pagenumber}{\normalfont\normalsize\sffamily}
\clearpairofpagestyles
\let\oldtitle\title
\renewcommand{\title}[1]{\oldtitle{#1}\newcommand{\theshorttitle}{#1}}
\newcommand{\shorttitle}[1]{\renewcommand{\theshorttitle}{#1}}
\let\oldauthor\author
\renewcommand{\author}[1]{\oldauthor{#1}\newcommand{\theshortauthor}{#1}}
\newcommand{\shortauthor}[1]{\renewcommand{\theshortauthor}{#1}}
\cohead{\xrfill[0.525ex]{0.6pt}~\theshorttitle~\xrfill[0.525ex]{0.6pt}}
\cehead{\xrfill[0.525ex]{0.6pt}~\theshortauthor~\xrfill[0.525ex]{0.6pt}}
\usepackage{natbib}
\usepackage{graphicx}
\usepackage[all]{xy}\UseComputerModernTips

\newcommand*{\arXiv}[1]{URL \bgroup\color{blue}\href{https://arxiv.org/abs/#1}{https://arxiv.org/abs/#1}\egroup}
\newcommand*{\doi}[1]{\bgroup\color{blue}\href{https://dx.doi.org/#1}{doi:#1}\egroup}
\newcommand*{\email}[1]{\bgroup\color{blue}\href{mailto:#1}{#1}\egroup}
\renewcommand*{\url}[1]{\bgroup\color{blue}\href{#1}{#1}\egroup}
\newcommand*{\ppara}[1]{\noindent\textbf{\textsf{#1}}\,\,}
\renewcommand{\cite}{\citep}

\newcommand*{\defeq}{:=}
\newcommand*{\Expect}{\mathbb{E}}

\newcommand*{\Naturals}{\mathbb{N}}
\newcommand*{\quark}{\setbox0\hbox{$x$}\hbox to\wd0{\hss$\cdot$\hss}}
\newcommand{\rd}{\textup{d}}
\newcommand*{\Reals}{\mathbb{R}}
\newcommand*{\Suldin}{Sul$'\!$din}
\newcommand{\todo}[1]{\bgroup\bfseries\color{red}#1\egroup}

\newcommand{\normal}{\mathcal{N}}
\newcommand*{\probson}[1]{\mathcal{P}_{#1}}

\newcommand*{\latentspace}{\mathcal{U}}
\newcommand*{\latentvar}{u}
\newcommand*{\dataspace}{\mathcal{Y}}
\newcommand*{\datavar}{y}
\newcommand*{\datamap}{Y}
\newcommand*{\qoispace}{\mathcal{Q}}
\newcommand*{\qoivar}{q}
\newcommand*{\qoimap}{Q}
\newcommand*{\nummap}{B}
\newcommand*{\pnummap}{\beta}

\newcommand{\arXivOmit}[1]{}
\newcommand{\STCOOmit}[1]{#1}

\newtheorem{theorem}{Theorem}[section]

\theoremstyle{definition}
\newtheorem{definition}[theorem]{Definition}
\newtheorem{example}[theorem]{Example}

\numberwithin{equation}{section}
\numberwithin{figure}{section}
\numberwithin{table}{section}

\begin{document}

\title{A Modern Retrospective on\\Probabilistic Numerics}
\shorttitle{A Modern Retrospective on Probabilistic Numerics}

\author{%
	C.\ J.\ Oates\footnote{Newcastle University, Herschel Building, Newcastle upon Tyne, NE1 7RU, UK, \email{chris.oates@ncl.ac.uk}} \footnote{Alan Turing Institute, British Library, 96 Euston Road, London NW1 2DB, UK, \email{coates@turing.ac.uk}}
	\and
	T.\ J.\ Sullivan\footnote{Institute of Mathematics, Freie Universit{\"a}t Berlin, Arnimallee 6, 14195 Berlin, Germany, \email{t.j.sullivan@fu-berlin.de}} \footnote{Zuse Institute Berlin, Takustra{\ss}e 7, 14195 Berlin, Germany, \email{sullivan@zib.de}}%
}
\shortauthor{C.\ J.\ Oates and T.\ J.\ Sullivan}

\date{\today}

\maketitle

\begin{abstract}
	\ppara{Abstract:}
	This article attempts to place the emergence of probabilistic numerics as a mathematical-statistical research field within its historical context and to explore how its gradual development can be related both to applications and to a modern formal treatment.
We highlight in particular the parallel contributions of {\Suldin} and Larkin in the 1960s and how their pioneering early ideas have reached a degree of maturity in the intervening period, mediated by paradigms such as average-case analysis and information-based complexity.
We provide a subjective assessment of the state of research in probabilistic numerics and highlight some difficulties to be addressed by future works.

	\smallskip
	
	\ppara{Keywords:}
	probabilistic numerics, scientific computation, reasoning under uncertainty, uncertainty quantification
	
	\smallskip
	
	\ppara{2010 Mathematics Subject Classification:}
	\mbox{62-03}, % STATISTICS > Historical
	\mbox{65-03}, % NUMERICAL ANALYSIS > Historical
	\mbox{01A60}, % History of mathematics and mathematicians > 20th century
	\mbox{01A65}, % History of mathematics and mathematicians > Contemporary
	\mbox{01A67} % History of mathematics and mathematicians > Future prospectives
\end{abstract}

\section{Introduction}
\label{intro}

The field of probabilistic numerics (PN), loosely speaking, attempts to provide a \emph{statistical} treatment of the errors and/or approximations that are made en route to the output of a deterministic numerical method, e.g.\ the approximation of an integral by quadrature, or the discretised solution of an ordinary or partial differential equation.
This decade has seen a surge of activity in this field.
In comparison with historical developments that can be traced back over more than a hundred years, the most recent developments are particularly interesting because they have been characterised by simultaneous input from multiple scientific disciplines:
mathematics, statistics, machine learning, and computer science.
The field has, therefore, advanced on a broad front, with contributions ranging from the building of over-arching general theory to practical implementations in specific problems of interest.
Over the same period of time, and because of increased interaction among researchers coming from different communities, the extent to which these developments were --- or were not --- presaged by twentieth-century researchers has also come to be better appreciated.

Thus, the time appears to be ripe for an update of the 2014 \emph{T{\"u}bingen Manifesto} on probabilistic numerics \cite{Hennig2014,Osborne2014d,Osborne2014c,Osborne2014b,Osborne2014a} and the position paper \cite{Hennig2015a} to take account of the developments between 2014 and 2019, an improved awareness of the history of this field, and a clearer sense of its future directions and potential.

In this article, we aim to summarise some of the history of probabilistic perspectives on numerics (Section~\ref{sec:history}), to place more recent developments into context (Section~\ref{sec:recent}), and to articulate a vision for future research in, and use of, probabilistic numerics (Section~\ref{sec:to_do_list}).

The authors are grateful to the participants of \emph{Prob Num 2018}, 11--13 April 2018, at the Alan Turing Institute, UK --- and in particular the panel discussants Oksana Chkrebtii, Philipp Hennig, Youssef Marzouk, Mike Osborne, and Houman Owhadi --- for many stimulating discussions on these topics.
However, except where otherwise indicated, the views that we present here are our own, and if we have misquoted or misrepresented the views of others, then the fault is entirely ours.

\section{Historical Developments}
\label{sec:history}

The first aim of this article is to reflect on the gradual emergence of probabilistic numerics as a research field.
The account in this section is not intended to be comprehensive in terms of the literature that is cited.
Rather, our aim is to provide an account of how the philosophical status of probabilistic approaches to numerical tasks has evolved, and in particular to highlight the parallel, pioneering, but often-overlooked contributions of {\Suldin} in the USSR and Larkin in the UK and Canada.

\subsection{Prehistory (--1959)}
\label{sec:prehistory}

The origins of PN can be traced to a discussion of probabilistic approaches to polynomial interpolation by Poincar\'{e} in his \emph{Calcul des Probabilit\'{e}s} (\cite[Ch.~21]{Poincare1896} and \cite[Ch.~25]{Poincare1912}).
Poincar\'{e} considered what, in modern terms, would be a particular case of a Gaussian infinite product measure prior on a function $f$, expressing it as a power series
\[
	f(x) = \sum_{k = 0}^{\infty} A_{k} x^{k}
\]
with independent normally-distributed coefficients $A_{k}$;
one is then given $n$ pointwise observations of the values of $f$ and seeks the probable values of $f(x)$ for another (not yet observed) value of $x$.

\begin{quotation}
	\noindent\textit{``Je suppose que l'on sache a priori que la fonction $f(x)$ est d{\'e}veloppable, dans une certain domaine, suivant les puissances croissantes des $x$,
	\[
		f(x) = A_{0} + A_{1} x + \dots .
	\]
	Nous ne savons rien sur les $A$, sauf que la probabilit\'{e} pour que l'un d'eux, $A_{i}$, soit compris entre certaines limites, $y$ et $y + \rd y$, est
	\[
		\sqrt{\frac{h_{i}}{\pi}} e^{- h_{i} y^{2}} \, \rd y .
	\]
	Nous connaissons par $n$ observations
	\begin{align*}
		f(a_{1}) & = B_{1} , \\
		f(a_{2}) & = B_{2} , \\
		\cdots \cdots & \cdots \cdots \\
		f(a_{n}) & = B_{n} .
	\end{align*}
	Nous cherchons la valeur probable de $f(x)$ pour une autre valeur de $x$.''}
	\cite[p.~292]{Poincare1912}
\end{quotation}

Note that, in using a Gaussian prior, Poincar\'{e} was departing from the Laplacian principle of indifference \cite{Laplace1812}, which would have mandated a uniform prior.\footnote{Indeed, while an improper uniform prior distribution on $\Reals$ makes sense for each $A_{k}$ individually, no such countably additive uniform measure (an ``infinite-dimensional Lebesgue measure'') can exist on $\Reals^{\infty}$ for $(A_{k})_{k = 0}^{\infty}$ \cite{Sudakov1959}.
That said, Poincar\'{e} does not impose any summability constraints on the $h_{i}$ either, so the covariance operator associated to his Gaussian prior may fail to be trace class.}

Poincar\'{e}'s analytical treatment predates the first digital multipurpose computers by decades, yet it clearly illustrates a non-trivial probabilistic perspective on a classic numerical task, namely function approximation by interpolation, a hybrid approach that is entirely in keeping with Poincar\'{e}'s reputation as one of the last universalist mathematicians \cite{Ginoux2013}.

However, our focus here is on the development of probabilistic numerical methods for use on a computer.
The limited nature of the earliest computers led authors to focus initially on the phenomenon of \emph{round-off error} \cite{Henrici1962,Hull1966,VonNeumannGoldstine1947}, whether of fixed-point or floating-point type, without any particular statistical \emph{inferential} motivation;
more recent contributions to the statistical study of round-off error include \cite{Barlow1985,Chatelin1990,Tienari1970}.
According to von Neumann and Goldstine, writing in 1947,
\begin{quotation}
	\noindent\textit{``[round-off errors] are strictly very complicated but uniquely defined number theoretical functions [of the inputs], yet our ignorance of their true nature is such that we best treat them as random variables.''}
	\cite[p.~1027]{VonNeumannGoldstine1947}.
\end{quotation}
Thus, von Neumann and Goldstine seem to have held a utilitarian view that probabilistic models in computation are useful shortcuts, simply easier to work with than the unwieldy deterministic truth.\footnote{Decades later, the discovery of chaotic dynamical systems would yield a similar conundrum:
after long enough time, one may as well assume that the system's state is randomly distributed according to its invariant measure, if it possesses one.}

Concerning the numerical solution of ordinary differential equations (ODEs), Henrici \cite{Henrici1962,Henrici1963} studied classical finite difference methods and derived expected values and covariance matrices for accumulated round-off error, under an assumption that individual round-off errors can be modelled as independent random variables.
In particular, given posited means and covariance matrices of the individual errors, Henrici demonstrated how these moments can be propagated through the computation of a finite difference method.
In contrast with more modern treatments, Henrici was concerned with the \emph{analysis} of an established numerical method and did not attempt to statistically \emph{motivate} the numerical method itself.

\subsection{The Parallel Contributions of Larkin and {\Suldin} (1959--1980)}
\label{subsec:LarkinSuldin}

One of the earliest attempts to motivate a numerical algorithm from a statistical perspective was due to Al$'\!$bert Valentinovich {\Suldin} (1924--1996), working at Kazan State University in the USSR (now Kazan Federal University in the Russian Federation) \cite{Norden1978,Zabotin1996}.
After first making contributions to the study of Lie algebras, towards the end of the 1950s {\Suldin} turned his attention to computational and applied mathematics, and in particular to probabilistic and statistical methodology.
His work in this direction led to the establishment of the Faculty of Computational Mathematics and Cybernetics (now Institute of Computational Mathematics and Information Technologies) in Kazan, of which he was the founding Dean.

\begin{figure}
	\centering
	\arXivOmit{\includegraphics[width = \linewidth]{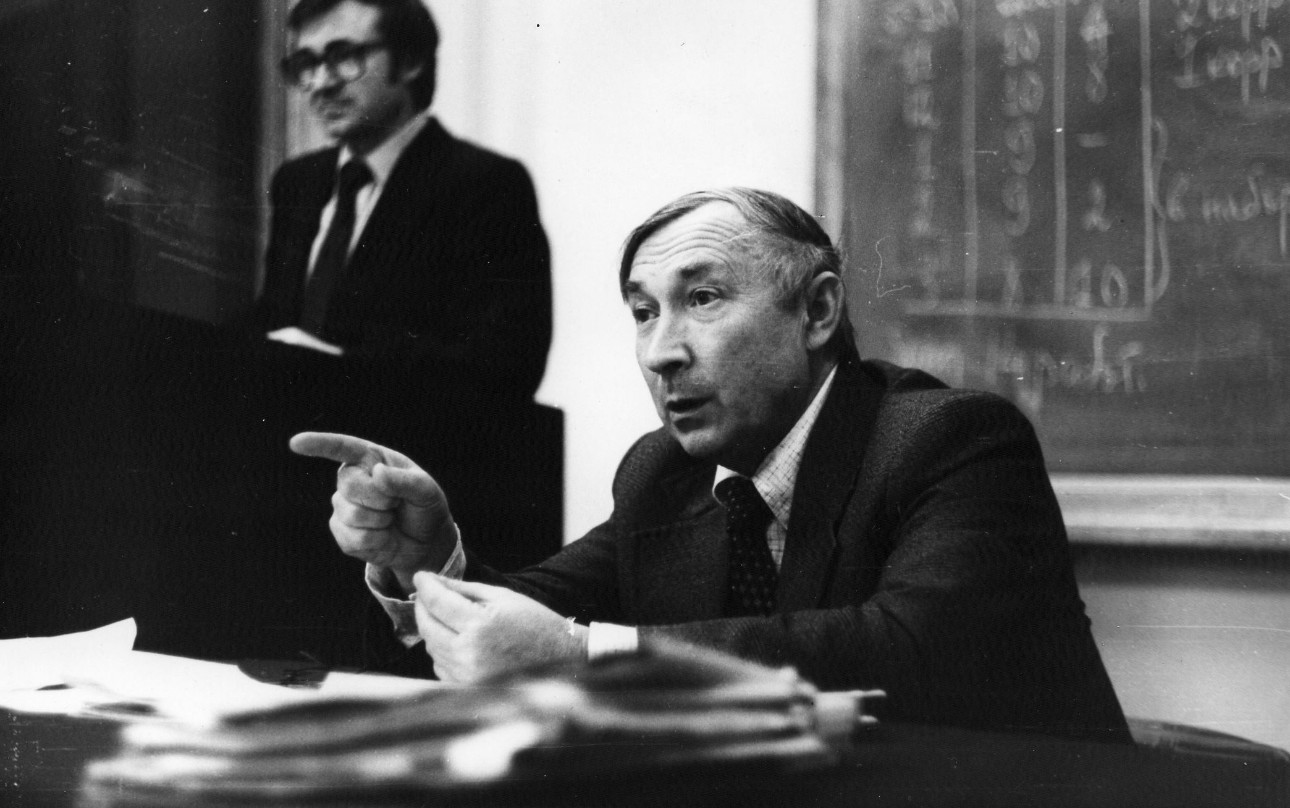}}
	\STCOOmit{\includegraphics[width = 0.5\linewidth]{figures/suldin2}}
	\caption{Al$'\!$bert Valentinovich {\Suldin} (1924--1996)
	\cite[reproduced with permission]{SuldinPicture}.}
	\label{fig:Suldin}
\end{figure}

{\Suldin} began by considering the problem of quadrature.
Suppose that we wish to approximate the definite integral $\int_{a}^{b} \latentvar(t) \, \rd t$ of a function $\latentvar \in \latentspace \defeq C^{0}([a, b]; \Reals)$, the space of continuous real-valued functions on $[a, b]$, under a statistical assumption that $(\latentvar(t) - \latentvar(a))_{t \in [a, b]}$ follows a standard Brownian motion (Wiener measure, $\mu_{\textup{W}}$).
For this task we receive pointwise data about the integrand $\latentvar$ in the form of the values of $\latentvar$ at $J \in \Naturals$ arbitrarily located nodes $t_{1}, \dots, t_{J} \in [a, b]$, although for convenience we assume that
\[
	a = t_{1} < t_{2} < \dots < t_{J} = b.
\]
In more statistical language, anticipating the terminology of Section~\ref{sec:Bayes_v_nonBayes}, our \emph{observed data} or \emph{information} concerning the integrand $\latentvar$ is $\datavar \defeq (t_{j}, \latentvar(t_{j}))_{j = 1}^{J}$, which takes values in the space $\dataspace \defeq ([a, b] \times \Reals)^{J}$.

Since $\mu_{\textup{W}}$ is a Gaussian measure and both the integral and pointwise evaluations of $\latentvar$ are linear functions of $\latentvar$, {\Suldin} \cite{Suldin1959,Suldin1960,Suldin1963b} showed by direct calculation that the quadrature rule $\nummap \colon \dataspace \to \Reals$ that minimises the mean squared error
\begin{equation}
	\label{eq:MSE}
	%\Expect_{\latentvar} \left[ \left| \int_{a}^{b} \latentvar(t) \, \rd t - \nummap \bigl( (t_{j}, \latentvar(t_{j}))_{j = 1}^{J} \bigr) \right|^{2} \right]
	\int_{\latentspace} \left| \int_{a}^{b} \latentvar(t) \, \rd t - \nummap \bigl( (t_{j}, \latentvar(t_{j}))_{j = 1}^{J} \bigr) \right|^{2} \, \mu_{\textup{W}} (\rd u)
\end{equation}
is the classical trapezoidal rule\footnote{Note that formulation \eqref{eq:trapezoidal_rule} of $\nummap_{\textup{tr}}$ emphasises the trapezoidal geometry being used to approximate the integral, whereas formulation \eqref{eq:trapezoidal_rule_bis} emphasises that the integrand need only be evaluated $J$ and not $2 J - 2$ times.}
\begin{align}
	\arXivOmit{ \notag
	& } \nummap_{\textup{tr}} \bigl( (t_{j}, z_{j})_{j = 1}^{J} \bigr) \arXivOmit{ \\ }
	\label{eq:trapezoidal_rule}
	& \arXivOmit{ \quad } \defeq \frac{1}{2} \sum_{j = 1}^{J - 1} (z_{j + 1} + z_{j}) (t_{j + 1} - t_{j}) \\
	\label{eq:trapezoidal_rule_bis}
	& \arXivOmit{ \quad } = z_{1} \frac{t_{2} - t_{1}}{2} + \sum_{j = 2}^{J - 1} z_{j} \frac{t_{j + 1} - t_{j - 1}}{2} + z_{J} \frac{t_{J} - t_{J - 1}}{2} ,
\end{align}
%%\textcolor{red}{Is the final line really more intuitive than (2), do you think? For me the final line is actually less intuitive than (2), but I am just a sample of one!}
%% Explained in a footnote.  it's to do with computational cost.
i.e.\ the definite integral of the piecewise linear interpolant of the observed data.
This result was a precursor to a sub-field of numerical analysis that became known as \emph{average-case analysis};
see Section~\ref{subsec:1980_2000}.

{\Suldin} was aware of the connection between his methods and statistical regression \cite{Suldin1963a} and conditional probability \cite{Suldin1963c}, although it is difficult to know whether he considered his work to be an expression of \emph{statistical inference} as such.
Indeed, since {\Suldin}'s methods were grounded in Hil\-bert space theory \cite{Suldin1968,Suldin1969}, the underlying mathematics (the linear conditioning of Gaussian measures on Hil\-bert spaces) is linear algebra which can be motivated without recourse to a probabilistic framework.

%%\textcolor{red}{I replaced ``Bayesian statistics'' with ``statistical inference'', since the trapezoidal rule is also the MLE, BLUE, etc, which are all estimators that exist in the frequentist world.}
%% OK!

In any case, {\Suldin}'s contributions were something entirely novel.
Up to this point, the role of statistics in numerical analysis was limited to providing insight into the \emph{performance} of a traditional numerical method.
The 1960s brought forth a new perspective, namely the stat\-istic\-ally-motivated \emph{design} of numerical methods.
Indeed,
\begin{quotation}
	\noindent\textit{``A.V.\ {\Suldin}'s 1969 habilitation thesis concern\-ed the development of probabilistic methods for the solution of problems in computational mathematics.
	His synthesis of two branches of mathematics turned out to be quite fruitful, and deep connections were discovered between the robustness of approximation formulae and their precision.
	Building on the general concept of an enveloping Hil\-bert space, A.V.\ {\Suldin} proved a projection theorem that enabled the solution of a number of approximation-theoretic problems.} \cite{Zabotin1996}
\end{quotation}

%{\Suldin}'s results on approximation quality were in the same spirit as the qualitative density results of Dunham Jackson \cite{Jackson1911}:
%
%\begin{theorem}
%	\label{thm:Jackson}
%	For each $r \in \Naturals$, there exists a universal constant $C(r)$ such that, for any $r$ times continuously differentiable $2 \pi$-periodic function $f$ with
%	\[
%		\max_{0 \leq k \leq r} \max_{0 \leq x \leq 2 \pi} | f^{(k)}(x) | \leq 1
%	\]
%	and any $n \in \Naturals$, there is a trigonometric polynomial $T_{n}$ of degree at most $n$ with
%	\[
%		\max_{0 \leq x \leq 2 \pi} | f(x) - T_{n}(x) | \leq \frac{C(r) \varpi(1/n, f^{(r)})}{n^{r}} ,
%	\]
%	where $\varpi(\delta, g)$ denotes the modulus of continuity of $g$ over a horizon $\delta > 0$.
%\end{theorem}

However, {\Suldin} was not alone in arriving at this point of view.
On the other side of the Iron Curtain, between 1957 and 1969, Frederick Michael (``Mike'') Larkin (1936--1982) worked for the UK Atomic Energy Authority in its laboratories at Harwell and Culham (the latter as part of the Computing and Applied Mathematics Group), as well as working for two years at Rolls Royce, England.
Following a parallel path to that of {\Suldin}, over the next decade Larkin would further blend numerical analysis and statistical thinking \cite{KuelbsLarkinWilliamson1972,Larkin1969,Larkin1972,Larkin1974,Larkin1979a,Larkin1979,Larkin1979b}, arguably laying the foundations on which PN would be developed.
At Culham, Larkin worked on building some of the first graphical calculators, called GHOST (short for \underline{g}raphical \underline{ou}tput sy\underline{st}em), and the GHOUL (\underline{g}raphical \underline{ou}tput \underline{l}anguage).
It can be speculated that an intimate familiarity with the computational limitations of GHOST and GHOUL may have motivated Larkin to seek a richer description of the numerical error associated to their output.

\begin{figure}
	\centering
	\arXivOmit{\includegraphics[width = \linewidth]{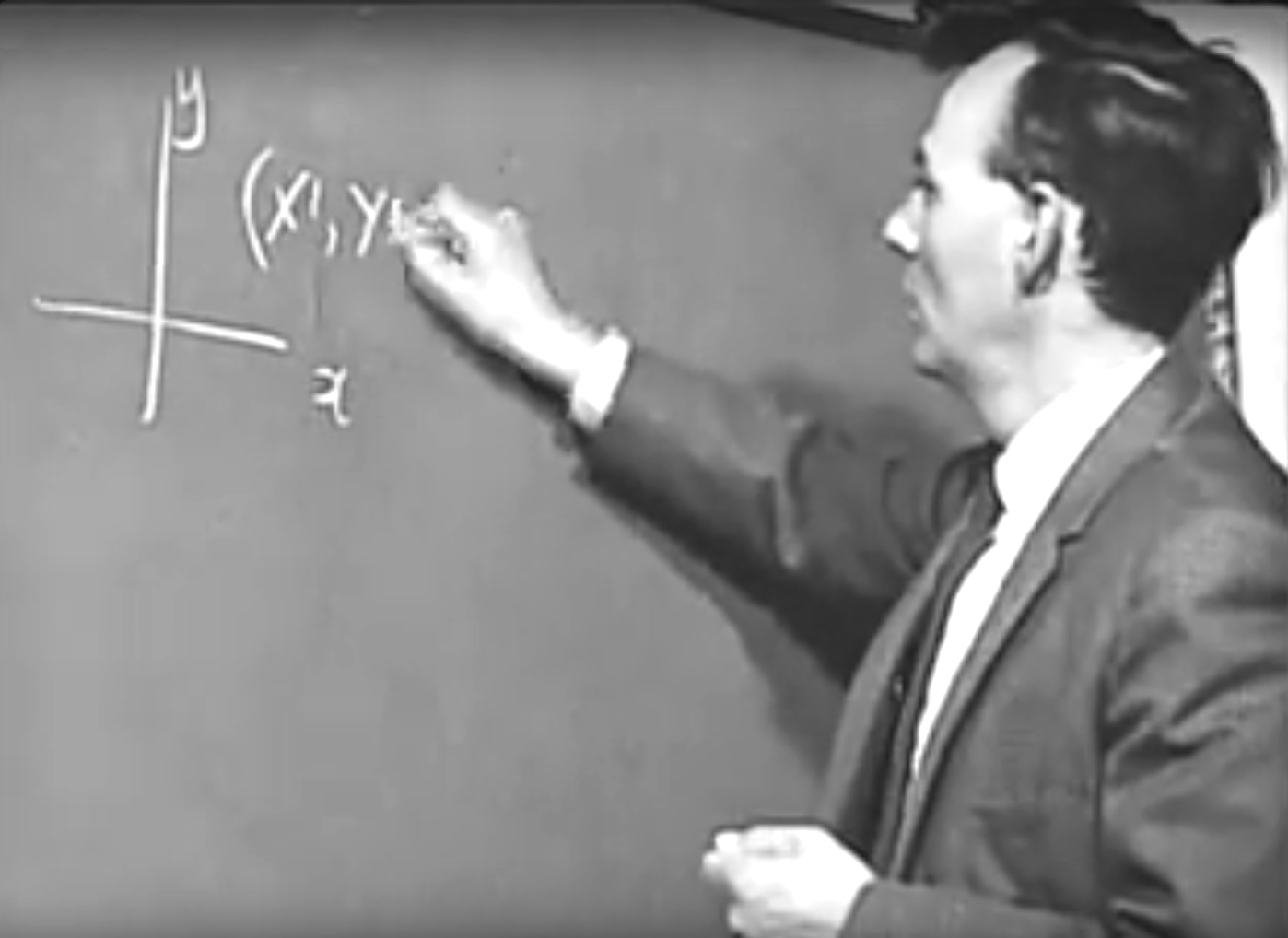}}
	\STCOOmit{\includegraphics[width = 0.5\linewidth]{figures/fml1}}
	\caption{Frederick Michael Larkin (1936--1982) \cite[reproduced with permission]{LarkinVideo}.}
	\label{fig:Larkin}
\end{figure}

The perspective developed by Larkin was fundamentally statistical and, in modern terminology, the probabilistic numerical methods he developed would be described as \emph{Bayesian}\footnote{Larkin used the term \emph{relative likelihood} for what we would recognise as a Bayesian prior \cite[Section 3.3]{Larkin1972}.
We may speculate, but cannot be sure, that such terminological differences are largely accidents of history.
Larkin was educated and did his early work exactly when the frequentist paradigm was starting to lose its dominance and Bayesian methods were starting to come back into fashion, driven by Cox's logical justification of the Bayesian paradigm \cite{Cox1946,Cox1961} and the development of theory, hardware, and software for methods like Markov chain Monte Carlo.
See \cite{Dale1999} for a comprehensive history of this area of statistics.},
which we discuss further in Section~\ref{sec:Bayes_v_nonBayes}.
Nevertheless, the pioneering nature of this research motivated Larkin to focus on specific numerical tasks, as opposed to establishing a unified framework.
In particular, he considered in detail the problems of approximating a non-negative function \cite{Larkin1969}, quadrature \cite{Larkin1972,Larkin1974}, and estimating the zeros of a complex function \cite{Larkin1979a,Larkin1979}.
In the context of the earlier numerical integration example of {\Suldin}, the alternative proposal of Larkin was to consider the Wiener measure as a prior, the information $(t_{j}, \latentvar(t_{j}))_{j=1}^J$ as (noiseless) data, and to output the posterior marginal for the integral $\int_{a}^{b} \latentvar(t) \, \rd t$.
That is, Larkin took the fundamental step of considering a distribution over the solution space of the numerical task to be the output of a computation --- this is what we would now recognise as the defining property of a \emph{probabilistic numerical method}:
\begin{quotation}
	\noindent\textit{``Among other things, this permits, at least in principle, the derivation of joint probability density functions for [both observed and unobserved] functionals on the space and also allows us to evaluate confidence limits on the estimate of a required functional (in terms of given values of other functionals).''% without any extra information about the norm of the function in question.
	} \cite{Larkin1972}\footnote{In this passage ``the estimate'' refers to the posterior mean in a linear-Gaussian set-up and ``confidence limit'' refers to what we would now call a highest-posterior-density credible interval.
	We suspect that the cultural dominance of frequentist statistics, in which estimators are reported alongside confidence intervals, led Larkin to adopt a similar presentation of the posterior --- though we emphasise that Larkin was fundamentally providing a Bayesian treatment.}
	%% \textcolor{blue}{should we abbreviate the ``without any extra information about the norm of the function in question'' bit? We don't comment on this part of the quote because it is a bit tangential.}
	%% Yep
\end{quotation}
Thus, in contrast to {\Suldin}'s description of the trapezoidal rule $\nummap_{\textup{tr}}$ from \eqref{eq:trapezoidal_rule} as a frequentist point estimator obtained from minimising \eqref{eq:MSE}, which just happens to produce an unbiased estimator with variance $\frac{1}{12} \sum_{j = 1}^{J - 1} (t_{j + 1} - t_{j})^{3}$, the Larkin viewpoint is to see the normal distribution 
\begin{equation}
	\label{eq:normal_trapezoidal}
	%\pnummap_{\textup{tr}} \bigl( (t_{j}, z_{j})_{j = 1}^{J} \bigr) \defeq 
	\normal \Biggl( \nummap_{\textup{tr}} \bigl( (t_{j}, z_{j})_{j = 1}^{J} \bigr) , \frac{1}{12} \sum_{j = 1}^{J - 1} (t_{j + 1} - t_{j})^{3} \Biggr)
\end{equation}
on $\Reals$ as the measure-valued output of a probabilistic quadrature rule, of which $\nummap_{\textup{tr}} \bigl( (t_{j}, z_{j})_{j = 1}^{J}$ is a convenient point summary.
Note also that the technical development in this pioneering work made fundamental contributions to the study of Gaussian measures on Hilbert spaces \cite{KuelbsLarkinWilliamson1972,Larkin1972}.

Larkin moved to Canada in 1969 to start work as a Consultant in Numerical Methods and Applied Mathematics within the Computing Centre and, subsequently in 1974, as Associate Professor in the Department of Computing and Information Science (now the School of Computing) at Queen's University in Kingston, Ontario.
He received tenure in 1977 and was promoted to full professor in 1980.
\begin{quotation}
	\noindent\textit{``He worked in isolation at Queen's in that few graduate students and fewer faculty members were aware of the nature of his research contributions to the field.
	[\dots]
	Michael pioneered the idea of using a probabilistic approach to give an alternative local approximation technique.
	In some cases this leads to the classical methods, but in many others leads to new algorithms that appear to have practical advantages over more classical methods.
	This work has finally begun to attract attention and I expect that the importance of his contribution will grow in time.''} \cite{Queens1982}
\end{quotation}

From our perspective, writing in 2019, it seems that {\Suldin} and Larkin were working in parallel but were ahead of their time.
Their probabilistic perspectives on approximation theory were similar, but limited to a Gaussian measure context.
Naturally, given the linguistic barriers and nearly disjoint publication cultures of their time, it would not have been easy for Larkin and {\Suldin} to be conversant with each other's work, though these barriers were not always as great as is sometimes thought \cite{Hollings2016}.
At least by 1972 \cite{Larkin1972}, Larkin was aware of and cited {\Suldin}'s work on minimal variance estimators for the values of linear functionals on Wiener space \cite{Suldin1959,Suldin1960}, but apparently did not know of {\Suldin}'s 1969 habilitation thesis, which laid out a broader agenda for the role of probability in numerics.
Conversely, Soviet authors writing in 1978 were aware of {\Suldin}'s influence on e.g.\ Ulf Grenander and Walter Freiberger at Brown University, but make no mention of Larkin \cite{Norden1978}.
{\Suldin}, for his part, at least as judged by his publication record, seems to have turned his attention to topics such as industrial mathematics (perhaps an ``easier sell'' in the production-oriented USSR \cite{Hollings2016}), mathematical biology, and of course the pressing concerns of faculty administration.

Finally, concerning the practicality of {\Suldin} and Larkin's ideas, one has to bear in mind the limited computational resources available at even cutting-edge facilities in the 1960s:\footnote{To first approximation, a single modern laptop has a hundred times the computing power of all five then-cutting-edge IBM System/360 Model 75J mainframe computers used for the ground support of the Apollo missions \cite{Manber2012}.}
probabilistic numerics was an idea ahead of its time, and the computational power needed to make it a reality simply did not exist.

\subsection{Optimal Numerical Methods are Bayes Rules (1980--1990)}
\label{subsec:1980_2000}

In the main, research contributions until 1990 continued to focus on deriving insight into traditional numerical methods through probabilistic analyses.
In particular, the \emph{average-case analysis} (ACA) of numerical methods received interest and built on the work of Kolmogorov \cite{Kolmogorov1936} and Sard \cite{Sard1963}.
In ACA the performance of a numerical method is assessed in terms of its \emph{average error} over an ensemble of numerical problems, with the ensemble being represented by a probability measure over the problem set;
a prime example is univariate quadrature with the average quadratic loss \eqref{eq:MSE} given earlier.
Root-finding, optimisation, etc.\ can all be considered similarly, and we defer to e.g.\ \cite{Ritter2000,Traub1983} for comprehensive treatments of this broad topic.

A traditional (deterministic) numerical method can also be regarded as a decision rule and the probability measure used in ACA can be used to instantiate the Bayesian decision-theoretic framework \cite{Berger1985}.
The average error is then recognised as the \emph{expected loss}, also called the \emph{risk}.
The fact that ACA is mathematically equivalent to Bayesian decision theory (albeit limited to the case of an experiment that produces a deterministic dataset) was noted in \cite{Kimeldorf1970a,Kimeldorf1970b,Parzen1970} and also in \cite{Larkin1970}.

Armed with an optimality criterion for a numerical method, it is natural to ask about the existence and performance of method(s) that minimise it.
Such methods are called \emph{average-case optimal} in ACA and are recognised as \emph{Bayes rules} or \emph{Bayes acts} in the decision-theoretic context.
A key result in this area is the insight of \cite{Kadane1985} that ACA-optimal methods coincide with (non-randomised) Bayes rules when the measure used to define the average error is the Bayesian prior;
for a further discussion of the relationships among these optimality criteria, including the Bayesian probabilistic numerical methods of Section~\ref{sec:Bayes_v_nonBayes}, see \cite{Cockayne_et_al_foundations,Oates_RICAM}.

Many numerical methods come in parametric families, being parametrised by e.g.\ the number of quadrature nodes, a mesh size, or a convergence tolerance.
For any ``sensible'' method, the error can be driven to zero by sending the parameter to infinity or zero as appropriate.
If one is prepared to pay an infinite computational cost, then essentially any method can be optimal!
Thus, when asking about the optimality of a numerical method, it is natural to consider the optimality of methods of a given computational cost or complexity.

With such concerns in mind, the field of \emph{information-based complexity} (IBC) \cite{Novak1988,Traub1983,Traub1980} developed simultaneously with ACA, with the aim of relating the computational complexity and optimality properties of algorithms to the available information on the unknowns, e.g.\ the partial nature of the information and any associated observational costs and errors.
%Despite the different names, there is in fact a great deal of overlap in terminology, concepts, and individual researchers between ACA and IBC.
For example, Smale \cite[Theorem~D]{Smale1985} compared the accuracies (with respect to mean absolute error) for a given cost of the Riemann sum, trapezoidal, and Simpson quadrature rules\footnote{On page 95 of the same paper, Smale highlighed Larkin's \cite{Larkin1972} as an \textit{``important earlier paper in this area''}.};
in the same paper, Smale also considered root-finding, optimisation via linear programming, and the solution of systems of linear equations.

The example of Bayesian quadrature was again discussed in detail by Diaconis \cite{Diaconis1988}, who repeated {\Suldin}'s observation that the posterior mean for $\int_{a}^{b} \latentvar(t) \, \rd t$ under the Wiener measure prior is the trapezoidal method \eqref{eq:trapezoidal_rule}, which is an ACA-optimal numerical method.
However, Diaconis posed a further question:
can other classical numerical integration methods, or numerical methods for other tasks, be similarly recovered as Bayes rules in a decision-theoretic framework?
For linear cubature methods, a positive and constructive answer was recently provided in \cite{Karvonen2018}, but the question remains open in general.

\subsection{Probabilistic Numerical Methods (1991--2009)}
\label{subsec:1991_2010}

After a period in which probabilistic numerical methods were all but forgotten, research interest was again triggered by contributions from \cite{Minka:2000,OHagan1991,RasmussenGhahramani:2003} on numerical integration, each to a greater or lesser extent a rediscovery of earlier work due to Larkin \cite{Larkin1972}.
In each case the output of computation was considered to be a probability distribution over the quantity of interest.

The 1990s saw an expansion in the PN agenda, first with early work on an area that was to become \emph{Bayesian optimisation} \cite{Mockus1974,Mockus1977,Mockus1989} and then with an entirely novel contribution on the numerical solution of ODEs by Skilling \cite{Skilling1992}.
Skilling presented a Bayesian\footnote{To be pedantic, the method of \cite{Skilling1992} does not satisfy the definition of a Bayesian PNM as given in Section~\ref{sec:Bayes_v_nonBayes}.
However, the method can be motivated as exact Bayesian inference under an approximate likelihood; see \cite{Wang:2018}.} perspective on the numerical solution of initial value problems of the form
\begin{align}
	\label{eq:IVP}
	\latentvar'(t) \equiv \frac{\rd \latentvar}{\rd t} & = f(t, \latentvar(t)) & & t \in [0, T] ,\\
	\notag
	\latentvar(0) & = \latentvar_{0} ,
\end{align}
and considered, for example, how regularity assumptions on $f$ should be reflected in correlation functions and the hypothesis space, how to choose a prior and likelihood, and potential sampling strategies.
Despite this work's then-new explicit emphasis on its Bayesian statistical character, Skilling himself considered his contributions to be quite natural:
\begin{quotation}
	\noindent\textit{``This paper arose from long exposure to Laplace/\linebreak[1]Cox/Jaynes probabilistic reasoning, combined with the University of Cambridge's desire that the author teach some (traditional) numerical analysis.
	The rest is common sense.
	[\dots]
	Simply, Bayes\-ian ideas are `in the air'.''} \cite{Skilling1992}
\end{quotation}

\subsection{Modern Perspective (2010--)}

The last two decades have seen an explosion of interest in \emph{uncertainty quantification} (UQ) for complex systems, with a great deal of research taking place in this area at the meeting point of applied mathematics, statistics, computational science, and application domains \cite{LeMaitreKnio2010,Smith:2014,Sullivan:2015}:
\begin{quotation}
	\noindent\textit{``UQ studies all sources of error and uncertainty, including the following:
	systematic and stochastic measurement error;
	ignorance;
	limitations of theoretical models;
	limitations of numerical representations of those models;
	limitations of the accuracy and reliability of computations, approximations, and algorithms;
	and human error.
	A more precise definition is UQ is the end-to-end study of the reliability of scientific inferences.''}
	\cite[p.~135]{GrandChallenges2009}
\end{quotation}
Since 2010, perhaps stimulated by this activity in the UQ community, a perspective on PN has emerged that sees PN part of UQ (broadly understood) and should be performed with a view to propagating uncertainty in computational pipelines.
This is discussed further in Sections~\ref{sec:hard_v_soft} and \ref{sec:Bayes_v_nonBayes}.

A notable feature of PN research since 2010 is the way that it has advanced on a broad front.
The topic of quadrature/cubature, in the tradition of {\Suldin} and Larkin, continues to be well represented:
see, e.g., \cite{Briol2018,Gunter2014,Karvonen2018,Oates2016,Osborne2012a,Osborne2012,Saerkkae2015,Xi2018} and \cite{Ehler2019,Jagadeeswaran2019,Karvonen2019b,Karvonen2019a}\arXivOmit{ in this special issue}.
The Bayesian approach to global optimisation continues to be widely used \cite{Chen2018,Snoek2012}, whilst probabilistic perspectives on quasi-Newton methods \cite{Hennig2013} and line search methods \cite{Mahsereci2015} have been put forward.
In the context of numerical linear algebra, \cite{Bartels2016,Cockayne2018,Hennig2015} and \cite{Bartels2019}\arXivOmit{ in this special issue} have approached the solution of a large linear system of equations as a statistical learning task and developed probabilistic alternatives to the classical conjugate gradient method.

Research has been particularly active in the development and analysis of statistical methods for the solution of ordinary and partial differential equations (ODEs and PDEs).
One line of research has sought to cast the solution of ODEs in the context of Bayesian filtering theory by building a Gaussian process (GP) regression model for the solution $\latentvar$ of the initial value problem of the form \eqref{eq:IVP}.
The observational data consists of the evaluations of the vector field $f$, interpreted as imperfect observations of the true time derivative $\latentvar'$, since one evaluates $f$ at the ``wrong'' points in space.
In this context, the key result is the Bayesian optimality of evaluating $f$ according to the classical Runge--Kutta (RK) scheme, so that the RK methods can be seen as point estimators of GP filtering schemes \cite{Kersting2016,Schober:2014wt,Schober:2016uh} and \cite{Tronarp2019}\arXivOmit{ in this special issue}.
Related iterative probabilistic numerical methods for ODEs include \cite{Abdulle2018,Chkrebtii2016,Conrad2017,KerstingSullivanHennig2018,Teymur2018,Teymur2016}.
The increased participation of mathematicians in the field has led to correspondingly deeper local and global convergence analysis of these methods in the sense of conventional numerical analysis, as in \cite{Conrad2017,KerstingSullivanHennig2018,Schober:2016uh,Teymur2018} and \cite{LieStuartSullivan2019}\arXivOmit{ in this special issue};
statistical principles for time step adaptivity have also been discussed, e.g.\ \cite{Chkrebtii2019}\arXivOmit{ in this special issue}.

For PDEs, resent research includes \cite{Chkrebtii2016,Cockayne2016,Cockayne2016_MaxEnt,Owhadi2015}, with these contributions making substantial use of reproducing kernel Hil\-bert space (RKHS) structure and Gaussian processes.
Unsurprisingly, given the deep connections between linear algebra and numerical methods for PDEs, the probabilistically-motivated theory of \emph{gamblets} for PDEs \cite{Owhadi2017sirev,Owhadi2017general,Owhadi2017jcp} has gone hand-in-hand with the development of fast solvers for structured matrix inversion and approximation problems \cite{Schaefer2017}; see also \cite{Yoo2019}\arXivOmit{ in this special issue}.

Returning to the point made at the beginning of this section, however, motivation for the development of probabilistic numerical methods has become closely linked to the traditional motivations of UQ (e.g.\  accurate and honest estimation of parameters of a so-called \emph{forward model}), with a role for PN due to the need to employ numerical methods to simulate from a forward model.
The idea to substitute a probability distribution in place of the (in general erroneous) output of a traditional numerical method can be used to prevent undue bias and over-confidence in the UQ task and is analogous to \emph{robust likelihood} methods in statistics \cite{Bissiri2016,Greco2008}.
This motivation is already present in \cite{Conrad2017}, and forms a major theme in \cite{Cockayne_et_al_foundations,Oates2017}.
Analysis of the impact of probabilistic numerical methods in simulation of the forward model within the context of Bayesian inversion has been provided in \cite{LieSullivanTeckentrup2018,StuartTeckentrup2018}.

\subsection{Related Fields and Their Development}

The field of PN did not emerge in isolation and the research cited above was undoubtedly influenced by parallel developments in mathematical statistics, some of which are now discussed.

First, the mathematical theory of \emph{optimal approximation using splines} was applied by Schoenberg \cite{Schoenberg1965,Schoenberg1966} and Karlin \cite{Karlin1969,Karlin1971,Karlin1972,Karlin1976} in the late 1960s and early 1970s to the linear problem of quadrature.
Larkin was aware of the work of Karlin, citing \cite{Karlin1969} in \cite{Larkin1974}.
However, the works cited above were not concerned with randomness and equivalent probabilistic interpretations were not discussed;
in contrast, the Bayesian interpretation of spline approximation was highlighted by \cite{Kimeldorf1970a}.

Second, the \emph{experimental design} literature of the late 1960s and early 1970s, including a sequence of contributions from Sacks and Ylvisacker \cite{SacksYlvisaker1968,SacksYlvisaker1970a,SacksYlvisaker1970b,SacksYlvisaker1966}, 
considered optimal selection of a design $0 \leq t_{1} < t_{2} < \dots < t_{J} \leq 1$ to minimise the covariance of the best linear estimator of $\beta$ given discrete observations of stochastic process
\begin{equation*}
	Y(t) = \sum_{i = 1}^{m} \beta_{i} \phi_{i} (t) + Z(t) ,
\end{equation*}
where $Z$ is a stochastic process with $\Expect [ Z(t) ] = 0$ and $\Expect [ Z(t)^{2} ] < \infty$, based on the data $\{ ( t_{j}, Y(t_{j}) ) \}_{j=1}^J$.
As such, the mathematical content of these works concerns optimal approximation in RKHSs, e.g.\ \cite[p.~2064, Theorem~1]{SacksYlvisaker1970a};
we note that Larkin \cite{Larkin1970} simultaneously considered optimal approximation in RKHSs.
However, the extent to which probability enters these works is limited to the measurement error process $Z$ that is entertained.

Third, the literature on \emph{emulation of black-box functions} that emerged in the late 1970s and 1980s, with contributions including \cite{OHagan1978,SacksWelchMitchellWynn1989}, provided Bayesian and frequentist statistical perspectives (respectively) on interpolation of a black-box function based on a finite number of function evaluations.
This literature did not present interpolation as an exemplar of other more challenging numerical tasks, such as the solution of differential equations, which could be similarly addressed but rather focused on the specific problem of black-box interpolation in and of itself.
The authors of \cite{SacksWelchMitchellWynn1989} were aware of the work of {\Suldin} but Larkin's work was not cited.
The challenges of proposing a suitable stochastic process model for a deterministic function were raised in the accompanying discussion of \cite{SacksWelchMitchellWynn1989} and were further discussed in \cite{Currin1991}.

\subsection{Conceptual Evolution --- A Summary}

To conclude and summarise this section, we perceive the following evolution of the concepts used in, and interpretation applied to, probability in numerical analysis:
\arXivOmit{\begin{enumerate}}
\STCOOmit{\begin{enumerate}[nosep]}
	\item In the traditional setting of numerical analysis, as seen circa 1950, all objects and operations are seen as being strictly deterministic.
	Even at that time, however, it was accepted by some that these deterministic objects are sometimes exceedingly complicated, to the extent that they may be treated as being stochastic, \`a la von Neumann and Goldstine.
	\item Sard and {\Suldin} considered the questions of optimal performance of a numerical method in, respectively, the worst-case and the average-case context.
	Though it is a fact that some of the average-case performance measures amount to variances of point estimators, they were not \emph{viewed} as such and in the early 1960s these probabilistic aspects were not a motivating factor.
	\item Larkin's innovation, in the late 1960s and early 1970s, was to formulate numerical tasks in terms of a joint distribution over latent quantities and quantities of interest, so that the quantity-of-interest output can be seen as a stochastic object.
	However, perhaps due to the then-prevailing statistical culture, Larkin summarised his posterior distributions using a point estimator accompanied by a credible interval.
	\item The fully modern viewpoint, circa 2019, is to explicitly think of the output as a probability measure to be realised, sampled, and possibly summarised.
\end{enumerate}

\section{Probabilistic Numerical Methods Come into Focus}
\label{sec:recent}

In this section we wish to emphasise how some of the recent developments mentioned in the previous section have brought greater clarity to the philosophical status of probabilistic numerics, clearing up some old points of disagreement or providing some standardised frameworks for the comparison of tasks and methods.

\subsection{A Means to an End, or an End in Themselves?}
\label{sec:hard_v_soft}

One aspect that has become clearer over the last few years, stimulated to some extent by disagreements between statisticians and numerical analysts over the role of probability in numerics, is that there are (at least) two distinct use cases or paradigms:
\arXivOmit{\begin{itemize}}
\STCOOmit{\begin{itemize}[nosep]}
	\item (P1) a probability-based \emph{analysis} of the performance of a (possibly classical) numerical method;
	\item (P2) a numerical method whose output carries the \emph{formal semantics} of some statistical inferential para\-digm (e.g.\ the Bayesian paradigm; cf.\ Section~\ref{sec:Bayes_v_nonBayes}).
\end{itemize}

Representatives of the first class of methods include \cite{Abdulle2018,Conrad2017}, which consider stochastic perturbations to explicit numerical integrators for ODEs in order to generate an ensemble of plausible trajectories for the unknown solution of the ODE.
In some sense, this can be viewed as a probabilistic sensitivity/stability analysis of a classical numerical method.
This first paradigm is also, clearly, closely related to ACA.

The second class of methods is exemplified by the Bayesian probabilistic numerical methods, discussed in \cite{Cockayne_et_al_foundations} and Section~\ref{sec:Bayes_v_nonBayes}.
We can further enlarge the second class to include those methods that only \emph{approximately} carry the appropriate semantics, e.g.\ because they are only approximately Bayesian, or only Bayesian for a particular quantity of interest or up to a finite time horizon, e.g.\ the filtering-based solvers for ODEs \cite{Kersting2016,KerstingSullivanHennig2018,Schober:2014wt,Schober:2016uh}.

Note that the second class of methods can also be pragmatically motivated, in the sense that formal statistical semantics enable techniques such as ANOVA to be brought to bear on the design and optimisation of a computational pipeline (to target the aspect of the computation that contributes most to uncertainty in the computational output) \cite{Hennig2015a}.
In this respect, statistical techniques can in principle supplement the expertise that is typically provided by a numerical analyst.

We note that paradigm (P1), with its close relationship to the longer-established field of ACA, tends to be more palatable to the classical numerical analysis community.
The typical, rather than worst-case, performance of a numerical method is of obvious practical interest \cite{Trefethen2008}.
Statisticians, especially practitioners of Bayesian and fiducial inference, are habitually more comfortable with paradigm (P2) than numerical analysts are.
As we remark in Section~\ref{sec:bridging}, this difference stems in part from a difference of opinion in which quantities are / can be regarded as ``random'' by the two communities;
this difference of opinion affects (P2) much more strongly than (P1).

\subsection{Bayesian Probabilistic Numerical Methods}
\label{sec:Bayes_v_nonBayes}

A recent research direction, which provides formal foundations for the approach pioneered by Larkin, is to interpret both traditional numerical methods and probabilistic numerical methods as particular solutions to an \emph{ill-posed inverse problem} \cite{Cockayne_et_al_foundations}.
Given that the latent quantities involved in numerical tasks are frequently functions, this development is in accordance with recent years' interest in non-parametric inversion in infinite-dimensional function spaces \cite{Stuart:2010,Sullivan:2015}.

From the point of view of \cite{Cockayne_et_al_foundations}, which echoes IBC, the common structure of numerical tasks such as quadrature, optimisation, and the solution of an ODE or PDE, is the following:
\arXivOmit{\begin{itemize}}
\STCOOmit{\begin{itemize}[nosep]}
	\item two known spaces:
	$\latentspace$, where the unknown latent variable lives,
	and $\qoispace$, where the quantity of interest lives;
	\item and a known function $\qoimap \colon \latentspace \to \qoispace$, a \emph{quantity-of-interest} function;
\end{itemize}
and the traditional role of the numerical analyst is to select/design
\arXivOmit{\begin{itemize}}
\STCOOmit{\begin{itemize}[nosep]}
	\item a space $\dataspace$, where data about the latent variable live;
	\item and two functions:
	$\datamap \colon \latentspace \to \dataspace$, an \emph{information operator} that acts on the latent variable to yield information,
	and $\nummap \colon \dataspace \to \qoispace$ such that $\nummap \circ \datamap \approx \qoimap$ in some sense to be determined.
\end{itemize}
With respect to this final point, Larkin \cite{Larkin1970} observed that there are many senses in which $\nummap \circ \datamap \approx \qoimap$.
One might ask, as Gaussian quadrature does, that the residual operator $R \defeq \nummap \circ \datamap - \qoimap$ vanish on a large enough finite-dimensional subspace of $\latentspace$;
one might ask, as worst-case analysis does, that $R$ be small in the supremum norm \cite{Sard1949};
one might ask, as ACA does, that $R$ be small in some integral norm against a probability measure on $\latentspace$.
In the chosen sense, numerical methods aim to make the following diagram approximately commute\footnote{Recall that a diagram such as \eqref{eq:ClassicalNM} or \eqref{eq:BPNM} is called \emph{commutative} if all routes that follow the arrows (functions) from any starting point to any endpoint yield the same result.
Thus, commutativity of \eqref{eq:ClassicalNM} means exactly that $\nummap(\datamap(\latentvar)) = \qoimap(\latentvar)$ for all $\latentvar \in \latentspace$.}:
\begin{equation}
	\label{eq:ClassicalNM}
	\xymatrix{
		\latentspace \ar[rr]^{\datamap} \ar[drr]_{\qoimap} & & \dataspace \ar@{-->}[d]^{\nummap} \\
		& & \qoispace
	}
\end{equation}

A statistician might say that a deterministic numerical method $\nummap \colon \dataspace \to \latentspace$ as described above uses observed data $\datavar \defeq \datamap(\latentvar)$ to give a \emph{point estimator} $\nummap(\datavar) \in \qoispace$ for a quantity of interest $\qoimap(\latentvar) \in \qoispace$ derived from a latent variable $\latentvar \in \latentspace$.

%\begin{figure}
%	\begin{center}
%		\xymatrix{
%			\latentspace \ar[rr]^{A} \ar[dr]_{Q} & & \dataspace \ar@{-->}[dl]^{b} \\
%			& \qoispace &
%		}
%	\end{center}
%	\caption{A deterministic numerical method $b \colon \dataspace \to \latentspace$ gives a point estimate $b(y) \in \qoispace$ for a quantity of interest $Q(u) \in \qoispace$ derived from a latent variable $u \in \latentspace$;
%	the observed data is $y = A(u)$.}
%	\label{fig:DNM}
%\end{figure}

\begin{example}
	\label{ex:integration_example}
	The general structure is exemplified by univariate quadrature, in which $\latentspace \defeq C^{0}([a, b]; \Reals)$, the information operator
	\[
		\datamap(\latentvar) \defeq (t_{j}, \latentvar(t_{j}))_{j = 1}^{J} \in \dataspace \defeq ([a, b] \times \Reals)^{J},
	\]
	corresponds to pointwise evaluation of the integrand at $J$ given nodes $a \leq t_{1} < \dots < t_{J} \leq b$, and the quantity of interest is
	\[
		\qoimap(\latentvar) \defeq \int_{a}^{b} \latentvar(t) \, \rd t \in \qoispace \defeq \Reals.
	\]
	Thus, we are interested in the definite integral of $\latentvar$, and we estimate it using only the information $\datamap(\latentvar)$, which does not completely specify $\latentvar$.
	Notice that \emph{some but not all} quadrature methods $\nummap \colon \dataspace \to \qoispace$ construct an estimate of $\latentvar$ and then exactly integrate this estimate;
	Gaussian quadrature does this by polynomially interpolating the observed data $\datamap(\latentvar)$;
	by way of constrast, vanilla Monte Carlo builds no such functional estimate of $\latentvar$, since its estimate for the quantity of interest,
	\begin{equation}
		\label{eq:vanilla_Monte_Carlo}
		\nummap_{\textup{MC}} \bigl( (t_{j}, z_{j})_{j = 1}^{J} \bigr) = \frac{1}{J} \sum_{j = 1}^{J} z_{j} ,
	\end{equation}
	forgets the locations $t_{j}$ at which the integrand $\latentvar$ was evaluated and uses only the values $z_{j} \defeq \latentvar(t_{j})$ of $\latentvar$.
	(Of course, the \emph{accuracy} of $\nummap_{\textup{MC}}$ is based on the assumption that the nodes $t_{j}$ are uniformly distributed in $[a, b]$.)
\end{example}

This formal framework enables a precise definition of a probabilistic numerical method (PNM) to be stated \cite[Section~2]{Cockayne_et_al_foundations}.
Assume that $\latentspace$, $\dataspace$, and $\qoispace$ are measurable spaces, that $\datamap$ and $\qoimap$ are measurable maps, and let $\probson{\latentspace}$ etc.\ denote the corresponding sets of probability distributions on these spaces.
Let $\qoimap_{\sharp} \colon \probson{\latentspace} \to \probson{\qoispace}$ denote the push-forward\footnote{I.e.\ $\qoimap_{\sharp}\mu(S) = \mu (\qoimap^{-1} (S))$ for all measurable $S \subseteq \qoispace$} of the map $\qoimap$, and define $\datamap_{\sharp}$ etc.\ similarly.

\begin{definition}
	\label{defn:PNM}
	A \emph{probabilistic numerical method} for the estimation of a quantity of interest $\qoimap$ consists of an information operator $\datamap \colon \latentspace \to \dataspace$ and a map $\pnummap \colon \probson{\latentspace} \times \dataspace \to \probson{\qoispace}$, the latter being termed a \emph{belief update operator}.
\end{definition}

That is, given a belief $\mu$ about $\latentvar$, $\pnummap(\mu, \quark)$ converts observed data $\datavar \in \dataspace$ about $\latentvar$ into a belief $\pnummap(\mu, \datavar) \in \probson{\qoispace}$ about $\qoimap(\latentvar)$, as illustrated by the dashed arrow in the following (not necessarily commutative) diagram:
\begin{equation}
	\label{eq:PNM}
	\xymatrix{
		\probson{\latentspace} \ar[rr]^{\datamap_{\sharp}} \ar[drr]_{\qoimap_{\sharp}} & & \probson{\dataspace} \ar[d]^{\nummap_{\sharp}} & & \dataspace \ar@{-->}[lld]^{\pnummap(\mu, \quark)} \ar@{..>}[d]^{\nummap} \\
		& & \probson{\qoispace} & & \qoispace \ar@{..>}[ll]^{\delta}
	}
\end{equation}

As shown by the dotted arrows in \eqref{eq:PNM}, this perspective is general enough to contain classical numerical methods $\nummap \colon \dataspace \to \qoispace$ as the special case $\pnummap(\mu, \datavar) = \delta_{\nummap(\datavar)}$, where $\delta_{\qoivar} \in \probson{\qoispace}$ is the unit Dirac measure at $\qoivar \in \qoispace$.

One desideratum for a PNM $\pnummap$ is that its point estimators (e.g.\ mean, median, or mode) should be closely related to standard deterministic numerical methods $\nummap$.
This aspect is present in works such as \cite{Schober:2014wt}, which considers probabilistic ODE solvers with Runge--Kutta schemes as their posterior means, and \cite{Cockayne2016,Cockayne2016_MaxEnt}, which consider PDE solvers with the symmetric collocation method as the posterior mean.
However, this aspect is by no means universally stressed.

A second, natural, desideratum for a PNM $\pnummap$ is that the spread (e.g.\ the variance) of the distributional output should provide a fair reflection of the accuracy to which the quantity of interest is being approximated.
In the statistics literature this amounts to a deside for credible intervals to be \emph{well calibrated} \cite{Robins2006}.
In particular, one might desire that the distribution $\pnummap$ contract to the true value of $\qoimap(\latentvar)$ at an appropriate rate as the data dimension (e.g.\ the number of quadrature nodes) is increased.\footnote{Here we abuse notation slightly:
strictly speaking, we should refer not to one PNM $\pnummap$ with input data $y$ of varying dimension but to a \emph{one-parameter family} of PNMs $\pnummap_{J}$ parametrised by the data dimension $J$.}

Diagram \eqref{eq:ClassicalNM}, when it commutes, characterises the ``ideal'' classical numerical method $\nummap$;
there is, as yet, no closed loop in diagram \eqref{eq:PNM} involving $\pnummap$, which we would need in order to describe an ``ideal'' PNM $\pnummap$.
This missing map in \eqref{eq:PNM} is intimately related to the notion of a \emph{Bayesian} PNM as defined by \cite{Cockayne_et_al_foundations}.

The key insight is that, given a prior belief expressed as a probability distribution $\mu \in \probson{\latentspace}$ and the information operator $\datamap \colon \latentspace \to \dataspace$, a Bayesian practitioner has a privileged map from $\dataspace$ into $\probson{\latentspace}$ to add to diagram \eqref{eq:PNM}, namely the conditioning operator that maps any possible value $\datavar \in \dataspace$ of the observed data to the corresponding conditional distribution $\mu^{\datavar} \in \probson{\latentspace}$ for $\latentvar$ given $\datavar$.
In this situation, in contrast to the freedom\footnote{The large and rapidly-growing canon of PNMs, only some of which are cited in this article, is strong evidence of just how great this freedom is!} enjoyed by the designer of an arbitrary PNM, a Bayesian has no choice in her/his belief $\pnummap(\mu, \datavar)$ about $\qoimap(\latentvar)$:
it must be nothing other than the image under $\qoimap$ of $\mu^{\datavar}$.

\begin{definition}
	\label{defn:BPNM}
	A probabilistic numerical method is said to be \emph{Bayesian} for $\mu \in \probson{\latentspace}$ if,
	\[
		\pnummap(\mu, \datavar) = \qoimap_{\sharp} \mu^{\datavar} \text{ for $\datamap_{\sharp} \mu$-almost all $\datavar \in \dataspace$ .}
	\]
	In this situation $\mu$ is called a \emph{prior} (for $\latentvar$) and $\pnummap(\mu, \datavar)$ a \emph{posterior} (for $\qoimap(\latentvar)$).
\end{definition}

In other words, being Bayesian means that the following diagram commutes:
\begin{equation}
	\label{eq:BPNM}
	\xymatrix{
		\probson{\latentspace} \ar[drr]_{\qoimap_{\sharp}} & & & & \dataspace \ar@{-->}[dll]^{\datavar \mapsto \pnummap(\mu, \datavar)} \ar[llll]_{\datavar \mapsto \mu^{\datavar}} \\
		& & \probson{\qoispace} & &
		}
\end{equation}
Note that Definition~\ref{defn:BPNM} does not insist that a Bayesian PNM actually calculates $\mu^{\datavar}$ and then computes the push-forward;
only that the output of the PNM is equal to $\qoimap_{\sharp} \mu^{\datavar}$.
Thus, whether or not a PNM is Bayesian is specific to the quantity of interest $\qoimap$.
Note also that a PNM $\pnummap(\mu, \quark)$ can be Bayesian for some priors $\mu$ yet be non-Bayesian for other choices of $\mu$;
for details see \cite[Sec.~5.2]{Cockayne_et_al_foundations}.

To be more formal for a moment, in Definition~\ref{defn:BPNM} the conditioning operation $\datavar \mapsto \mu^{\datavar}$ is interpreted in the sense of a \emph{disintegration}, as advocated by \cite{Chang1997}.
This level of technicality is needed in order to make rigorous sense of the operation of conditioning on the $\mu$-negligible event that $\datamap(\latentvar) = \datavar$.
Thus,
\arXivOmit{\begin{itemize}}
\STCOOmit{\begin{itemize}[nosep]}
	\item for each $\datavar \in \dataspace$, $\mu^{\datavar} \in \probson{\latentspace}$ is supported only on those values of $\latentvar$ compatible with the observation $\datamap(\latentvar) = \datavar$, i.e.\ $\mu^{\datavar} ( \{ \latentvar \in \latentspace \mid \datamap(\latentvar) \neq \datavar \} ) = 0$;
	\item for any measurable set $E \subseteq \latentspace$, $\datavar \mapsto \mu^{\datavar}(E)$ is a measurable function from $\dataspace$ into $[0, 1]$ satisfying the \emph{reconstruction property}, or \emph{law of total probability},
	\[
		\mu(E) = \int_{\dataspace} \mu^{\datavar}(E) \, (\datamap_{\sharp} \mu) (\rd \datavar) .
	\]
\end{itemize}
Under mild conditions\footnote{Sufficient conditions are, e.g., that $\latentspace$ be a complete and separable metric space with its Borel $\sigma$-algebra (so that every $\mu \in \probson{\latentspace}$ is a Radon measure) and that the $\sigma$-algebra on $\dataspace$ be countably generated and contain all singletons.} such a disintegration always exists, and is unique up to modification on $\datamap_{\sharp} \mu$-null sets.

Observe that the fundamental difference between ACA (i.e.\ the probabilistic assessment of classical numerical methods) and Bayesianity of PNMs is that the former concerns the commutativity of diagram \eqref{eq:ClassicalNM} in the average (i.e.\ the left-hand half of diagram \eqref{eq:PNM}), whereas the latter concerns the commutativity of diagram \eqref{eq:BPNM}.

The prime example of a Bayesian PNM is the following example of \emph{kernel quad\-rature}, due to \cite{Larkin1972}:

\begin{example}
	\label{eg:kernel_quadrature}
	Recall the set-up of Example~\ref{ex:integration_example}.
	Take a Gaussian distribution $\mu$ on $C^{0}([a, b]; \Reals)$, with mean function $m \colon [a, b] \to \Reals$ and covariance function $k \colon [a, b]^2 \to \Reals$.
	Then, given the data
	\[
		y = (t_{j}, z_{j})_{j = 1}^{J} \equiv (t_{j}, \latentvar(t_{j}))_{j = 1}^{J} ,
	\]
	the disintegration $\mu^{\datavar}$ is again a Gaussian on $C^{0}([a, b]; \Reals)$ with mean and covariance functions
	\begin{align}
		\label{eq:kernel_quadrature_post_mean_func}
		m^{\datavar} (t) & = m(t) + k_{T} (t)^\top k_{TT}^{-1} (z_{T} - m_{T}) , \\
		\label{eq:kernel_quadrature_post_cov_func}
		k^{\datavar} (t, t') & = k(t, t') - k_{T} (t)^{\top} k_{TT}^{-1} k_{T}(t') ,
	\end{align}
%	where $[k_{T} (t)]_j \defeq k(t, t_{j})$, $[k_{TT}]_{i,j} \defeq k(t_{i}, t_{j})$, $[z_{T}]_{j} \defeq z_{j} \equiv \latentvar(t_{j})$, and $[m_{T}]_j \defeq m(t_{j})$.
	where $k_{T} \colon [a, b] \to \Reals^{J}$, $k_{TT} \in \Reals^{J \times J}$, $z_{T} \in \Reals^{J}$, and $m_{T} \in \Reals^{J}$ are given by
	\begin{align*}
		[k_{T} (t)]_j & \defeq k(t, t_{j}) , &
		[k_{TT}]_{i,j} & \defeq k(t_{i}, t_{j}) , \\
		[z_{T}]_{j} & \defeq z_{j} \equiv \latentvar(t_{j}) , &
		[m_{T}]_j & \defeq m(t_{j}) .
	\end{align*}
	The Bayesian PNM output $\pnummap(\mu, \datavar)$, i.e.\ the push-forward $\qoimap_{\sharp} \mu^{\datavar}$, is a Gaussian on $\Reals$ with mean $\overline{m}^{y}$ and variance $(\overline{\sigma}^{y})^{2}$ given by integrating \eqref{eq:kernel_quadrature_post_mean_func} and \eqref{eq:kernel_quadrature_post_cov_func} respectively, i.e.
	\begin{align*}
		\overline{m}^{y} & = \int_{a}^{b} m(t) \, \rd t + \left[ \int_{a}^{b} k_{T} (t) \, \rd t \right]^\top k_{TT}^{-1} (z_{T} - m_{T}) , \\
		(\overline{\sigma}^{y})^{2} & = \int_{a}^{b} \int_{a}^{b} k(t, t') \, \rd t \, \rd t' \arXivOmit{ \\
		& \phantom{=} \quad } - \left[ \int_{a}^{b} k_{T} (t) \, \rd t \right]^\top k_{TT}^{-1} \left[\int_{a}^{b} k_{T} (t') \, \rd t' \right] .
	\end{align*}
	From a practical perspective, $k$ is typically taken to have a parametric form $k_\theta$ and the parameters $\theta$ are adjusted in a data-dependent manner, for example to maximise the marginal likelihood of the information $\datavar$ under the Gaussian model.
	
	One may also seek point sets that minimise the posterior variance $(\overline{\sigma}^{y})^{2}$ of the estimate of the integral.
	For the Brownian covariance kernel $k(t, t') = \min(t, t')$, the posterior $\qoimap_{\sharp} \mu = \normal(\overline{m}^{y}, (\overline{\sigma}^{y})^{2})$ for $\int_{a}^{b} \latentvar(t) \, \rd t$ is given by \eqref{eq:normal_trapezoidal}, the variance of which is clearly minimised by an equally-spaced point set $\{ t_{j} \}_{j=1}^{J}$.
	For more general kernels $k$, an early reference for selecting the point set $\{ t_{j} \}_{j=1}^{J}$ to minimise $(\overline{\sigma}^{y})^{2}$ is \cite{OHagan1991}.
\end{example}

This perspective, in which the Bayesian update is singled out from other possible belief updates, is reminiscent of foundational discussions such as \cite{Bissiri2016,Zellner:1988}.
Interestingly, about half of the papers published on PN can be viewed as being (at least approximately) Bayes\-ian; see the survey in the supplement of \cite{Cockayne_et_al_foundations}.
This includes the work of Larkin, though, as previously mentioned, Larkin himself did not use the terminology of the Bayesian framework.
Quite aside from questions of computational cost, non-Bayesian methods come into consideration because the requirement to be fully Bayes\-ian can impose non-trivial constraints on the design of a practical numerical method, particularly for problems with a causal aspect or ``time's arrow'';
this point was discussed in detail for the numerical solution of ODEs in \cite{Wang:2018}.

As well as providing a clear formal benchmark, \cite[Section~5]{Cockayne_et_al_foundations} argue that a key advantage of Bayesian probabilistic numerical methods is that they are \emph{closed under composition}, so that the output of a computational pipeline composed of Bayesian probabilistic numerical methods will inherit Bayesian semantics itself.
This is analogous to the Markov condition that underpins directed acyclic graphical models \cite{Lauritzen1996} and may be an advantageous property in the context of large and/or distributed computational codes --- an area where performing a classical numerical analysis can often be difficult.
For non-Bayesian PNMs it is unclear how these can/should be combined, but we note an analogous discussion of statistical ``models made of modules'' in the recent work of \cite{Jacob:2017} (who observe, like \cite{Owhadi2015_brittleness}, that strictly Bayesian models can be brittle under model misspecification, whereas non-Bayesianity confers additional robustness) and also the numerical analysis of probabilistic forward models in Bayesian inverse problems in \cite{LieSullivanTeckentrup2018}.

\section{Discussion and Outlook}
\label{sec:to_do_list}

\begin{quotation}
	\noindent\textit{``Det er vanskeligt at spaa, is{\ae}r naar det g{\ae}lder Fremtiden.''}
	[Danish proverb]
\end{quotation}

As it stands in 2019, our view is that there is much to be excited about.
An intermittent stream of ad hoc observations and proposals, which can be traced back to the pioneering work of Larkin and {\Suldin}, has been unified under the banner of probabilistic numerics \cite{Hennig2015a} and solid statistical foundations have now been established \cite{Cockayne_et_al_foundations}.
In this section we comment on some of the most important aspects of research that remain to be addressed.

\subsection{Killer Apps}

The most successful area of research to date has been on the development of Bayesian methods for global optimisation \cite{Snoek2012}, which have become standard to the point of being embedded into commercial software \cite{BayesOptMatlab} and deployed in realistic \cite{Acerbi2018,Paul2018} and indeed high-profile \cite{Chen2018} applications.
Other numerical tasks have yet to experience the same level of practical interest, though we note applications of probabilistic methods for cubature in computer graphics \cite{Marques2013} and tracking \cite{Pruher2018}, as well as applications of probabilistic numerical methods in medical tractography \cite{Hauberg2015} and nonlinear state estimation \cite{Oates2017} in an industrial context.

It has been suggested that probabilistic numerics is likely to experience the most success in addressing numerical tasks that are fundamentally difficult \cite{Owen2019}.
One area that we highlight, in particular, in this regard is the solution of high-dimensional PDEs.
There is considerable current interest in the deployment of neural networks as a substitute for more traditional numerical methods in this context, e.g.\ \cite{Sirignano2018}, and the absence of interpretable error indicators for neural networks is a strong motivation for the development of more formal probabilistic numerical methods for this task.
We note also that nonlinear PDEs in particular are prone to non-uniqueness of solutions.
For some problems, physical reasoning may be used to choose among the various solutions, from the probabilistic or statistical perspective lack of uniqueness presents no fundamental philosophical issues:
the multiple solutions are simply multiple maxima of a likelihood, and the prior is used to select among them, as in e.g.\ the treatment of Painlev\'{e}'s transcendents in \cite{Cockayne_et_al_foundations}.

It has also been noted that the probabilistic approach provides a promising paradigm for the analysis of rounding error in mixed-precision calculations, where classical bounds ``do not provide good estimates of the size of the error, and in particular [\dots] overestimate the error growth, that is, the asymptotic dependence of the error on the problem size'' \cite{Higham2018}.

\subsection{Adaptive Bayesian Methods}

The presentation of a PNM in Section~\ref{sec:Bayes_v_nonBayes} did not permit \emph{adaptation}.
It has been rigorously established that for \emph{linear} problems adaptive methods (e.g., in quadrature, sequential selection of the notes $t_{j}$) do not outperform non-adaptive methods according to certain performance metrics such as worst-case error \cite[Section~3.2]{Wozniakowski1985}.
However, adaptation is known to be advantageous in general for \emph{nonlinear} problems \cite[Section~3.8]{Wozniakowski1985}.
At a practical level, adaptation is usually an essential component in the development of stopping rules that enable a numerical method to terminate after an error indicator falls below a certain user-specified level.
An analysis of adaptive PNMs would constitute a non-trivial generalisation of the framework of \cite{Cockayne_et_al_foundations}, who limited attention to static directed acyclic graph representation of conditional dependence structure.
The generalisation to adaptive PNM necessitates the use of graphical models with a natural filtration, as exemplified by a dynamical Bayesian network \cite{Murphy2002}.

It has been suggested that numerical analysis is a natural use case for \emph{empirical Bayes methods} \cite{Carlin2000,Casella1985}, as opposed to related -- but usually more computationally intensive -- approaches such as hierarchical modelling and cross-validation.
Empirical Bayes methods can be characterised as a specific instance of adaptation in which the observed data are used not only for inference but also to form a point estimator for the prior.
For example, in a quadrature setting, the practitioner is in the fortunate position of being able to use evaluations of the integrand $\latentvar$ both to estimate the regularity of $\latentvar$ and the value of the integral.
Empirical Bayesian methods are explored in \cite{Schober:2016uh} and in \cite{Jagadeeswaran2019}\arXivOmit{ in this special issue}.

\subsection{Design of Probabilistic Numerical Methods}

Paradigmatic questions in the IBC literature are those of (i) an optimal information operator $\datamap$ for a given task, and (ii) the optimal numerical method $\nummap$ for a given task, given information of a known type \cite{Traub1983}.
In the statistical literature, there is also a long history of Bayesian optimal experimental design, in parametric and non-parametric contexts \cite{Lindley1956,Piiroinen2005}.
The extent to which these principles can be used to design optimal numerical methods automatically (rather than by inspired guesswork on the mathematician's part, \`a la Larkin) remains a major open question, analogous to the automation of statistical reasoning envisioned by Wald and subsequent commentators on his work \cite{OwhadiScovel2017}.

\subsection{Probabilistic Programming}

The theoretical foundations of probabilistic numerics have now been laid, but at present a library of compatible code has not been developed.
In part, this is due to the amount of work needed in order to make a numerical implementation reliable and efficient, and in this respect PN lies far behind classical numerical analysis at present.
Nevertheless, we anticipate that such efforts will be undertaken in coming years, and will lead to the wider adoption of probabilistic numerical methods.
In particular, we are excited at the prospect of integrating probabilistic numerical methods into a probabilistic programming language, e.g.\ \cite{Carpenter2017}, where tools from functional programming and category theory can be exploited in order to automatically compile codes built from probabilistic numerical methods \cite{Scibior2015}.

\subsection{Bridging the Numerics--Statistics Gap}
\label{sec:bridging}

\begin{quotation}
	\noindent\textit{``Numerical analysts and statisticians are both in the business of estimating parameter values from incomplete information.
	The two disciplines have separately developed their own approaches to formalizing strangely similar problems and their own solution techniques;
	the author believes they have much to offer each other.''} \cite{Larkin1979b}
\end{quotation}

A major challenge faced by researchers in this area is the interdisciplinary gap between numerical analysts on the one hand and statisticians on the other.
Though there are some counterexamples, as a first approximation it is true to say that classically-trained numerical analysts lack deep knowledge of probability or statistics, and classically-trained statisticians are not well versed in numerical topics such as convergence and stability analysis.
Indeed, not only do these two communities take interest in different questions, they often fail to even see the point of the other group's expertise and approaches to their common problems.

A caricature of this mutual incomprehension is the following:
A numerical analyst will quite rightly point out that almost all problems have numerical errors that are provably non-Gaussian, not least because s/he can exhibit a rigorous a-priori or a-posteriori error bound.
Therefore, to the numerical analyst it seems wholly inappropriate to resort to Gaussian models for any purpose at all;
these are often the statistician's first models of choice, though they should not be the last.
This non-paradox was explained in detail by Larkin in \cite{Larkin1974}.
(As a side note, it seems to us from our discussions that numerical analysts are happier to discuss the modelling of \emph{errors} than the \emph{latent quantities} which they regard as fixed, whereas statisticians seems to have the opposite preference;
this is a difference in views that echoes the famous frequentist--subjectivist split in statistics.)
The numerical analyst also wonders why, in the presence of an under-resolved integral, the practitioner does not simply apply an adaptive quadrature scheme and run it until an \textit{a posteriori} global error indicator falls below a pre-set tolerance.

We believe that these difficulties are not fundamental and can be overcome by a more careful statement of the approach being taken to address the numerical task.
In particular, the meeting ground for the numerical analysts and statisticians, and the critical arena of application for PN, consists of problems that \emph{cannot} be run to convergence more cheaply than quantifying the uncertainties of the coarse solution --- or, at least, where there is an interesting cost-v.-accuracy tradeoff to be had, which is a central enabling factor for multilevel methods \cite{Giles2015}.

More generally, we are encouraged to see that epistemic uncertainty is being used once again and an analytical device in numerical analysis in the sense originally described by von Neumann and Goldstine \cite{VonNeumannGoldstine1947};
see e.g.\ \cite{Higham2018}.

\subsection{Summary}
\label{subsec:summary}

The first aim of this article was to better understand probabilistic numerics through its historical development.
Aside from the pioneering work of Larkin, it was only in the 1990s that probabilistic numerical methods --- i.e.\ algorithms returning a probability distribution as their output --- were properly developed.
A unified vision of probabilistic computation was powerfully presented in \cite{Hennig2015a} and subsequently formalised in \cite{Cockayne_et_al_foundations}.

The second aim of this article was to draw a distinction between PN as a means to an end, as a form of probabilistic sensitivity / stability analysis, and PN as an end in itself.
In particular, we highlighted the Bayesian sub-class of PNMs as being closed under composition, a property that makes these particularly well suited for use in UQ;
we also remarked that many problems --- for reasons of problem structure, computational cost, or robustness to model misspecification --- call for methods that are not formally Bayesian.

Finally, we highlighted areas for further development, which we believe will be essential if the full potential of probabilistic numerics highlighted in \cite{Hennig2015a} is to be realised.
From our perspective, the coming to fruition of this vision will require demonstrable success on problems that were intractable with the computational resources of previous decades and a wider acceptance of Larkin's observation quoted above, with which we wholeheartedly agree:
numerical analysts and statisticians are indeed in the same business and do have much to offer one other!

\section*{Acknowledgements}
\addcontentsline{toc}{section}{Acknowledgements}

The authors wish to express their sincere thanks to Paul Constantine for highlighting the contribution of \cite{VonNeumannGoldstine1947};
to Matthew Larkin, Graham Larkin, Cherrilyn Yalin and Selim Akl for assisting with the historical content in Section~\ref{subsec:LarkinSuldin} and giving permission to reproduce Figure~\ref{fig:Larkin};
to Sergey Mosin and Kazan Federal University for permission to reproduce Figure~\ref{fig:Suldin};
to Milena Kremakova for help with translating \cite{Norden1978,Zabotin1996} from Russian into English and for assisting in communications with Kazan Federal University;
and to Ilse Ipsen and four anonymous reviewers for their thoughtful comments.

CJO was supported by the Lloyd's Register Foundation programme on Data-Centric Engineering at the Alan Turing Institute, London, UK.

TJS was supported by the Freie Universit{\"a}t Berlin within the Excellence Initiative of the German Research Foundation (DFG), and by the DFG Collaborative Research Centre 1114 ``Scaling Cascades in Complex Systems''.

% The following acknowledgment text should be included
This material was developed, in part, at the \textit{Prob Num 2018} workshop hosted by the Lloyd's Register Foundation programme on Data-Centric Engineering at the Alan Turing Institute, UK, and supported by the National Science Foundation, USA, under Grant DMS-1127914 to the Statistical and Applied Mathematical Sciences Institute.
Any opinions, findings, conclusions or recommendations expressed in this material are those of the author(s) and do not necessarily reflect the views of the above-named funding bodies and research institutions.

\bibliographystyle{abbrvnat}
\addcontentsline{toc}{section}{References}
\bibliography{references.bib}

\def\cprime{$'\!$}
\begin{thebibliography}{156}
\providecommand{\natexlab}[1]{#1}
\providecommand{\url}[1]{\texttt{#1}}
\expandafter\ifx\csname urlstyle\endcsname\relax
  \providecommand{\doi}[1]{doi: #1}\else
  \providecommand{\doi}{doi: \begingroup \urlstyle{rm}\Url}\fi

\bibitem[Abdulle and Garegnani(2018)]{Abdulle2018}
A.~Abdulle and G.~Garegnani.
\newblock Random time step probabilistic methods for uncertainty quantification
  in chaotic and geometric numerical integration, 2018.
\newblock \arXiv{1801.01340}.

\bibitem[Acerbi(2018)]{Acerbi2018}
L.~Acerbi.
\newblock Variational {Bayesian} {Monte} {Carlo}.
\newblock In \emph{32nd Conference on Neural Information Processing Systems
  (NeurIPS 2018)}, 2018.
\newblock URL
  \url{https://papers.nips.cc/paper/8043-variational-bayesian-monte-carlo}.

\bibitem[Barlow and Bareiss(1985)]{Barlow1985}
J.~L. Barlow and E.~H. Bareiss.
\newblock Probabilistic error analysis of {Gaussian} elimination in floating
  point and logarithmic arithmetic.
\newblock \emph{Computing}, 34\penalty0 (4):\penalty0 349--364, 1985.
\newblock URL \url{https://doi.org/10.1007/BF02251834}.

\bibitem[Bartels and Hennig(2016)]{Bartels2016}
S.~Bartels and P.~Hennig.
\newblock Probabilistic approximate least-squares.
\newblock In \emph{Proceedings of the 19th International Conference on
  Artificial Intelligence and Statistics}, volume~51 of \emph{Proceedings of
  Machine Learning Research}, pages 676--684, 2016.
\newblock URL \url{http://proceedings.mlr.press/v51/bartels16.pdf}.

\bibitem[Bartels et~al.(2019)Bartels, Cockayne, Ipsen, and Hennig]{Bartels2019}
S.~Bartels, J.~Cockayne, I.~C.~F. Ipsen, and P.~Hennig.
\newblock Probabilistic linear solvers: {A} unifying view.
\newblock \emph{Stat. Comp.}, 2019.
\newblock To appear. \arXiv{1810.03398}.

\bibitem[Berger(1985)]{Berger1985}
J.~O. Berger.
\newblock \emph{Statistical {Decision} {Theory} and {Bayesian} {Analysis}}.
\newblock Springer Series in Statistics. Springer-Verlag, New York, second
  edition, 1985.
\newblock URL \url{https://doi.org/10.1007/978-1-4757-4286-2}.

\bibitem[Bissiri et~al.(2016)Bissiri, Holmes, and Walker]{Bissiri2016}
P.~G. Bissiri, C.~C. Holmes, and S.~G. Walker.
\newblock A general framework for updating belief distributions.
\newblock \emph{J. Roy. Stat. Soc. Ser. B}, 78\penalty0 (5):\penalty0
  1103--1130, 2016.
\newblock URL \url{https://doi.org/10.1111/rssb.12158}.

\bibitem[Briol et~al.(2019)Briol, Oates, Girolami, Osborne, and
  Sejdinovic]{Briol2018}
F.-X. Briol, C.~J. Oates, M.~Girolami, M.~A. Osborne, and D.~Sejdinovic.
\newblock Probabilistic integration: A role in statistical computation? (with
  discussion and rejoinder).
\newblock \emph{Stat. Sci.}, 34\penalty0 (1):\penalty0 1--22, 2019.
\newblock URL \url{https://doi.org/10.1214/18-STS660}.

\bibitem[Carlin and Louis(2000)]{Carlin2000}
B.~P. Carlin and T.~A. Louis.
\newblock Empirical {Bayes}: past, present and future.
\newblock \emph{J. Amer. Stat. Assoc.}, 95\penalty0 (452):\penalty0 1286--1289,
  2000.
\newblock URL \url{https://doi.org/10.2307/2669771}.

\bibitem[Carpenter et~al.(2017)Carpenter, Gelman, Hoffman, Lee, Goodrich,
  Betancourt, Brubaker, Guo, Li, and Riddell]{Carpenter2017}
B.~Carpenter, A.~Gelman, M.~Hoffman, D.~Lee, B.~Goodrich, M.~Betancourt,
  M.~Brubaker, J.~Guo, P.~Li, and A.~Riddell.
\newblock Stan: {A} probabilistic programming language.
\newblock \emph{J. Stat. Software}, 76\penalty0 (1), 2017.
\newblock URL \url{https://doi.org/10.18637/jss.v076.i01}.

\bibitem[Casella(1985)]{Casella1985}
G.~Casella.
\newblock An introduction to empirical {Bayes} data analysis.
\newblock \emph{Amer. Stat.}, 39\penalty0 (2):\penalty0 83--87, 1985.
\newblock URL \url{https://doi.org/10.2307/2682801}.

\bibitem[Chang and Pollard(1997)]{Chang1997}
J.~T. Chang and D.~Pollard.
\newblock Conditioning as disintegration.
\newblock \emph{Stat. Neerlandica}, 51\penalty0 (3):\penalty0 287--317, 1997.
\newblock URL \url{https://doi.org/10.1111/1467-9574.00056}.

\bibitem[Chatelin and Brunet(1990)]{Chatelin1990}
F.~Chatelin and M.-C. Brunet.
\newblock A probabilistic round-off error propagation model. {Application} to
  the eigenvalue problem.
\newblock In \emph{Reliable numerical computation}, Oxford Sci. Publ., pages
  139--160. Oxford Univ. Press, New York, 1990.

\bibitem[Chen et~al.(2018)Chen, Huang, Wang, Antonoglou, Schrittwieser, Silver,
  and de~Freitas]{Chen2018}
Y.~Chen, A.~Huang, Z.~Wang, I.~Antonoglou, J.~Schrittwieser, D.~Silver, and
  N.~de~Freitas.
\newblock {Bayesian} optimization in {AlphaGo}, 2018.
\newblock \arXiv{1812.06855}.

\bibitem[Chkrebtii and Campbell(2019)]{Chkrebtii2019}
O.~A. Chkrebtii and D.~A. Campbell.
\newblock Adaptive step-size selection for state-space based probabilistic
  differential equation solvers.
\newblock \emph{Stat. Comp.}, 2019.
\newblock To appear.

\bibitem[Chkrebtii et~al.(2016)Chkrebtii, Campbell, Calderhead, and
  Girolami]{Chkrebtii2016}
O.~A. Chkrebtii, D.~A. Campbell, B.~Calderhead, and M.~A. Girolami.
\newblock Bayesian solution uncertainty quantification for differential
  equations.
\newblock \emph{Bayesian Anal.}, 11\penalty0 (4):\penalty0 1239--1267, 2016.
\newblock URL \url{https://doi.org/10.1214/16-BA1017}.

\bibitem[Cockayne et~al.(2016)Cockayne, Oates, Sullivan, and
  Girolami]{Cockayne2016}
J.~Cockayne, C.~Oates, T.~J. Sullivan, and M.~Girolami.
\newblock Probabilistic meshless methods for partial differential equations and
  {Bayesian} inverse problems, 2016.
\newblock \arXiv{1605.07811}.

\bibitem[Cockayne et~al.(2017)Cockayne, Oates, Sullivan, and
  Girolami]{Cockayne2016_MaxEnt}
J.~Cockayne, C.~Oates, T.~J. Sullivan, and M.~Girolami.
\newblock Probabilistic numerical methods for {PDE}-constrained {Bayesian}
  inverse problems.
\newblock In G.~Verdoolaege, editor, \emph{Proceedings of the
  36\textsuperscript{th} International Workshop on Bayesian Inference and
  Maximum Entropy Methods in Science and Engineering}, volume 1853 of \emph{AIP
  Conference Proceedings}, pages 060001--1--060001--8, 2017.
\newblock URL \url{https://doi.org/10.1063/1.4985359}.

\bibitem[Cockayne et~al.(2019{\natexlab{a}})Cockayne, Oates, Sullivan, and
  Girolami]{Cockayne_et_al_foundations}
J.~Cockayne, C.~Oates, T.~J. Sullivan, and M.~Girolami.
\newblock {Bayesian} probabilistic numerical methods.
\newblock \emph{SIAM Rev.}, 2019{\natexlab{a}}.
\newblock To appear. \arXiv{1702.03673}.

\bibitem[Cockayne et~al.(2019{\natexlab{b}})Cockayne, Oates, Ipsen, and
  Girolami]{Cockayne2018}
J.~Cockayne, C.~J. Oates, I.~Ipsen, and M.~Girolami.
\newblock A {Bayesian} conjugate gradient method.
\newblock \emph{Bayesian Anal.}, 2019{\natexlab{b}}.
\newblock To appear. \arXiv{1801.05242}.

\bibitem[Conrad et~al.(2017)Conrad, Girolami, S\"{a}rkk\"{a}, Stuart, and
  Zygalakis]{Conrad2017}
P.~R. Conrad, M.~Girolami, S.~S\"{a}rkk\"{a}, A.~Stuart, and K.~Zygalakis.
\newblock Statistical analysis of differential equations: introducing
  probability measures on numerical solutions.
\newblock \emph{Stat. Comp.}, 27\penalty0 (4):\penalty0 1065--1082, 2017.
\newblock URL \url{https://doi.org/10.1007/s11222-016-9671-0}.

\bibitem[Cox(1946)]{Cox1946}
R.~T. Cox.
\newblock Probability, frequency and reasonable expectation.
\newblock \emph{Amer. J. Phys.}, 14\penalty0 (1):\penalty0 1--13, 1946.
\newblock URL \url{https://doi.org/10.1119/1.1990764}.

\bibitem[Cox(1961)]{Cox1961}
R.~T. Cox.
\newblock \emph{The {Algebra} of {Probable} {Inference}}.
\newblock The Johns Hopkins Press, Baltimore, MD, 1961.

\bibitem[Currin et~al.(1991)Currin, Mitchell, Morris, and
  Ylvisaker]{Currin1991}
C.~Currin, T.~Mitchell, M.~Morris, and D.~Ylvisaker.
\newblock Bayesian prediction of deterministic functions, with applications to
  the design and analysis of computer experiments.
\newblock \emph{J. Amer. Stat. Assoc.}, 86\penalty0 (416):\penalty0 953--963,
  1991.
\newblock URL \url{https://doi.org/10.1080/01621459.1991.10475138}.

\bibitem[Dale(1999)]{Dale1999}
A.~I. Dale.
\newblock \emph{A {History} of {Inverse} {Probability}: {From} {Thomas} {Bayes}
  to {Karl} {Pearson}}.
\newblock Sources and Studies in the History of Mathematics and Physical
  Sciences. Springer-Verlag, New York, second edition, 1999.
\newblock URL \url{https://doi.org/10.1007/978-1-4419-8652-8}.

\bibitem[Diaconis(1988)]{Diaconis1988}
P.~Diaconis.
\newblock Bayesian numerical analysis.
\newblock In \emph{Statistical decision theory and related topics, {IV},
  {Vol.}\ 1 ({West} {Lafayette}, {Ind.}, 1986)}, pages 163--175, New York,
  1988. Springer.
\newblock URL \url{https://doi.org/10.1007/978-1-4613-8768-8\_20}.

\bibitem[Ehler et~al.(2019)Ehler, Gr\"af, and Oates]{Ehler2019}
M.~Ehler, M.~Gr\"af, and C.~J. Oates.
\newblock Optimal {Monte} {Carlo} integration on closed manifolds.
\newblock \emph{Stat. Comp.}, 2019.
\newblock To appear. \arXiv{1707.04723}.

\bibitem[Giles(2015)]{Giles2015}
M.~B. Giles.
\newblock Multilevel {Monte} {Carlo} methods.
\newblock \emph{Acta Numer.}, 24:\penalty0 259--328, 2015.
\newblock URL \url{https://doi.org/10.1017/S096249291500001X}.

\bibitem[Ginoux and Gerini(2013)]{Ginoux2013}
J.~M. Ginoux and C.~Gerini.
\newblock \emph{Henri {Poincar{\'e}}: {A} {Biography} {Through} the {Daily}
  {Papers}}.
\newblock World Scientific, 2013.
\newblock URL \url{https://doi.org/10.1142/8956}.

\bibitem[Greco et~al.(2008)Greco, Racugno, and Ventura]{Greco2008}
L.~Greco, W.~Racugno, and L.~Ventura.
\newblock Robust likelihood functions in {Bayesian} inference.
\newblock \emph{J. Stat. Plann. Inference}, 138\penalty0 (5):\penalty0
  1258--1270, 2008.
\newblock URL \url{https://doi.org/10.1016/j.jspi.2007.05.001}.

\bibitem[Gunter et~al.(2014)Gunter, Osborne, Garnett, Hennig, and
  Roberts]{Gunter2014}
T.~Gunter, M.~A. Osborne, R.~Garnett, P.~Hennig, and S.~J. Roberts.
\newblock Sampling for inference in probabilistic models with fast {Bayesian}
  quadrature.
\newblock In \emph{Advances in Neural Information Processing Systems 27}, pages
  2789--2797, 2014.
\newblock URL
  \url{https://papers.nips.cc/paper/5483-sampling-for-inference-in-probabilistic-models-with-fast-bayesian-quadrature}.

\bibitem[Hauberg et~al.(2015)Hauberg, Schober, Liptrot, Hennig, and
  Feragen]{Hauberg2015}
S.~Hauberg, M.~Schober, M.~Liptrot, P.~Hennig, and A.~Feragen.
\newblock A random {Riemannian} metric for probabilistic shortest-path
  tractography.
\newblock volume 9349 of \emph{Lecture Notes in Computer Science}, pages
  597--604. 2015.
\newblock URL \url{https://doi.org/10.1007/978-3-319-24553-9\_73}.

\bibitem[Hennig(2014)]{Hennig2014}
P.~Hennig.
\newblock Roundtable in {T}\"{u}bingen, 2014.
\newblock URL
  \url{http://www.probnum.org/2014/08/22/Roundtable-2014-in-Tuebingen/}.

\bibitem[Hennig(2015)]{Hennig2015}
P.~Hennig.
\newblock Probabilistic interpretation of linear solvers.
\newblock \emph{SIAM J. Optim.}, 25\penalty0 (1):\penalty0 234--260, 2015.
\newblock URL \url{https://doi.org/10.1137/140955501}.

\bibitem[Hennig and Kiefel(2013)]{Hennig2013}
P.~Hennig and M.~Kiefel.
\newblock Quasi-{Newton} methods: {A} new direction.
\newblock \emph{J. Mach. Learn. Research}, 14\penalty0 (Mar):\penalty0
  843--865, 2013.
\newblock URL
  \url{http://www.jmlr.org/papers/volume14/hennig13a/hennig13a.pdf}.

\bibitem[Hennig et~al.(2015)Hennig, Osborne, and Girolami]{Hennig2015a}
P.~Hennig, M.~A. Osborne, and M.~Girolami.
\newblock Probabilistic numerics and uncertainty in computations.
\newblock \emph{Proc. R. Soc. A}, 471\penalty0 (2179):\penalty0 20150142, 2015.
\newblock URL \url{https://doi.org/10.1098/rspa.2015.0142}.

\bibitem[Henrici(1962)]{Henrici1962}
P.~Henrici.
\newblock \emph{Discrete {V}ariable {M}ethods in {O}rdinary {D}ifferential
  {E}quations}.
\newblock John Wiley \& Sons, Inc., New York-London, 1962.

\bibitem[Henrici(1963)]{Henrici1963}
P.~Henrici.
\newblock \emph{Error {P}ropagation for {D}ifference {M}ethod}.
\newblock John Wiley \& Sons, Inc., New York-London, 1963.

\bibitem[Higham and Mary(2018)]{Higham2018}
N.~J. Higham and T.~Mary.
\newblock A new approach to probabilistic rounding error analysis.
\newblock Technical report, University of Manchester, 2018.
\newblock URL \url{http://eprints.maths.manchester.ac.uk/2673/1/paper.pdf}.

\bibitem[Hollings(2016)]{Hollings2016}
C.~D. Hollings.
\newblock \emph{Scientific {Communication} {Across} the {Iron} {Curtain}}.
\newblock SpringerBriefs in History of Science and Technology. Springer, Cham,
  2016.
\newblock URL \url{https://doi.org/10.1007/978-3-319-25346-6}.

\bibitem[Hull and Swenson(1966)]{Hull1966}
T.~E. Hull and J.~R. Swenson.
\newblock Tests of probabilistic models for the propagation of roundoff errors.
\newblock \emph{Comm. ACM}, 9:\penalty0 108--113, 1966.
\newblock URL \url{https://doi.org/10.1145/365170.365212}.

\bibitem[Jacob et~al.(2017)Jacob, Murray, Holmes, and Robert]{Jacob:2017}
P.~E. Jacob, L.~M. Murray, C.~C. Holmes, and C.~P. Robert.
\newblock Better together? {Statistical} learning in models made of modules,
  2017.
\newblock \arXiv{1708:08719}.

\bibitem[Jagadeeswaran and Hickernell(2018)]{Jagadeeswaran2019}
R.~Jagadeeswaran and F.~J. Hickernell.
\newblock Fast automatic {Bayesian} cubature using lattice sampling, 2018.
\newblock \arXiv{1809.09803}.

\bibitem[Kadane and Wasilkowski(1985)]{Kadane1985}
J.~B. Kadane and G.~W. Wasilkowski.
\newblock Average case {$\varepsilon$}-complexity in computer science. {A}
  {Bayesian} view.
\newblock In \emph{Bayesian Statistics, 2 ({V}alencia, 1983)}, pages 361--374.
  North-Holland, Amsterdam, 1985.

\bibitem[Karlin(1969)]{Karlin1969}
S.~Karlin.
\newblock Best quadrature formulas and interpolation by splines satisfying
  boundary conditions.
\newblock In \emph{Approximations with {Special} {Emphasis} on {Spline}
  {Functions} ({Proc.} {Sympos.} {Univ.} of {Wisconsin}, {Madison}, {Wis.},
  1969)}, pages 447--466. Academic Press, New York, 1969.

\bibitem[Karlin(1971)]{Karlin1971}
S.~Karlin.
\newblock Best quadrature formulas and splines.
\newblock \emph{J. Approx. Theory}, 4:\penalty0 59--90, 1971.
\newblock URL \url{https://doi.org/10.1016/0021-9045(71)90040-2}.

\bibitem[Karlin(1972)]{Karlin1972}
S.~Karlin.
\newblock On a class of best nonlinear approximation problems.
\newblock \emph{Bull. Amer. Math. Soc.}, 78:\penalty0 43--49, 1972.
\newblock URL \url{https://doi.org/10.1090/S0002-9904-1972-12842-8}.

\bibitem[Karlin(1976)]{Karlin1976}
S.~Karlin.
\newblock \emph{Studies in Spline Functions and Approximation Theory}, chapter
  On a class of best nonlinear approximation problems and extended monosplines,
  pages 19--66.
\newblock Academic Press, New York, 1976.

\bibitem[Karvonen et~al.(2018{\natexlab{a}})Karvonen, Kanagawa, and
  S\"{a}rkk\"{a}]{Karvonen2019b}
T.~Karvonen, M.~Kanagawa, and S.~S\"{a}rkk\"{a}.
\newblock On the positivity and magnitudes of {Bayesian} quadrature weights,
  2018{\natexlab{a}}.
\newblock \arXiv{1812.08509}.

\bibitem[Karvonen et~al.(2018{\natexlab{b}})Karvonen, Oates, and
  S\"{a}rkk\"{a}]{Karvonen2018}
T.~Karvonen, C.~J. Oates, and S.~S\"{a}rkk\"{a}.
\newblock A {Bayes}--{Sard} cubature method.
\newblock In \emph{32nd Conference on Neural Information Processing Systems
  (NeurIPS 2018)}, 2018{\natexlab{b}}.
\newblock URL
  \url{http://papers.nips.cc/paper/7829-a-bayes-sard-cubature-method}.

\bibitem[Karvonen et~al.(2019)Karvonen, S\"{a}rkk\"{a}, and
  Oates]{Karvonen2019a}
T.~Karvonen, S.~S\"{a}rkk\"{a}, and C.~J. Oates.
\newblock Symmetry exploits for {Bayesian} cubature methods.
\newblock \emph{Stat. Comp.}, 2019.
\newblock To appear. \arXiv{1809.09803}.

\bibitem[{Kazan Federal University}()]{SuldinPicture}
{Kazan Federal University}.
\newblock URL \url{https://kpfu.ru/portal/docs/F\_261937733/suldin2.jpg}.
\newblock Accessed December 2018.

\bibitem[Kersting and Hennig(2016)]{Kersting2016}
H.~Kersting and P.~Hennig.
\newblock Active uncertainty calibration in {Bayesian} {ODE} solvers.
\newblock In \emph{Proceedings of the 32nd Conference on Uncertainty in
  Artificial Intelligence (UAI 2016)}, pages 309--318, 2016.
\newblock URL \url{http://www.auai.org/uai2016/proceedings/papers/163.pdf}.

\bibitem[Kersting et~al.(2018)Kersting, Sullivan, and
  Hennig]{KerstingSullivanHennig2018}
H.~Kersting, T.~J. Sullivan, and P.~Hennig.
\newblock Convergence rates of {Gaussian} {ODE} filters, 2018.
\newblock \arXiv{1807.09737}.

\bibitem[Kimeldorf and Wahba(1970{\natexlab{a}})]{Kimeldorf1970a}
G.~S. Kimeldorf and G.~Wahba.
\newblock A correspondence between {Bayesian} estimation on stochastic
  processes and smoothing by splines.
\newblock \emph{Ann. Math. Stat.}, 41:\penalty0 495--502, 1970{\natexlab{a}}.
\newblock URL \url{https://doi.org/10.1214/aoms/1177697089}.

\bibitem[Kimeldorf and Wahba(1970{\natexlab{b}})]{Kimeldorf1970b}
G.~S. Kimeldorf and G.~Wahba.
\newblock Spline functions and stochastic processes.
\newblock \emph{Sankhy\={a} Ser. A}, 32:\penalty0 173--180, 1970{\natexlab{b}}.
\newblock URL \url{https://www.jstor.org/stable/25049652}.

\bibitem[Kolmogorov(1936)]{Kolmogorov1936}
A.~N. Kolmogorov.
\newblock {\"{U}ber} die beste {Ann\"{a}herung} von {Funktionen} einer
  gegebenen {Funktionenklasse}.
\newblock \emph{Ann. of Math. (2)}, 37\penalty0 (1):\penalty0 107--110, 1936.
\newblock URL \url{https://doi.org/10.2307/1968691}.

\bibitem[Kuelbs et~al.(1972)Kuelbs, Larkin, and
  Williamson]{KuelbsLarkinWilliamson1972}
J.~Kuelbs, F.~M. Larkin, and J.~A. Williamson.
\newblock Weak probability distributions on reproducing kernel {Hilbert}
  spaces.
\newblock \emph{Rocky Mountain J. Math.}, 2\penalty0 (3):\penalty0 369--378,
  1972.
\newblock URL \url{https://doi.org/10.1216/RMJ-1972-2-3-369}.

\bibitem[Laplace(1812)]{Laplace1812}
P.~S. Laplace.
\newblock \emph{Th{\'e}orie {Analytique} des {Probabilit{\'e}s}}.
\newblock Courcier, Paris, 1812.

\bibitem[Larkin(1969)]{Larkin1969}
F.~M. Larkin.
\newblock Estimation of a non-negative function.
\newblock \emph{BIT Num. Math.}, 9\penalty0 (1):\penalty0 30--52, 1969.
\newblock URL \url{https://doi.org/10.1007/BF01933537}.

\bibitem[Larkin(1970)]{Larkin1970}
F.~M. Larkin.
\newblock Optimal approximation in {Hilbert} spaces with reproducing kernel
  functions.
\newblock \emph{Math. Comp.}, 24:\penalty0 911--921, 1970.
\newblock URL \url{https://doi.org/10.2307/2004625}.

\bibitem[Larkin(1972)]{Larkin1972}
F.~M. Larkin.
\newblock Gaussian measure in {Hilbert} space and applications in numerical
  analysis.
\newblock \emph{Rocky Mountain J. Math.}, 2\penalty0 (3):\penalty0 379--421,
  1972.
\newblock URL \url{https://doi.org/10.1216/RMJ-1972-2-3-379}.

\bibitem[Larkin(1974)]{Larkin1974}
F.~M. Larkin.
\newblock Probabilistic error estimates in spline interpolation and quadrature.
\newblock In \emph{Information Processing 74 (Proc.\ IFIP Congress, Stockholm,
  1974)}, pages 605--609, Amsterdam, 1974. North-Holland.

\bibitem[Larkin(1979{\natexlab{a}})]{Larkin1979}
F.~M. Larkin.
\newblock A modification of the secant rule derived from a maximum likelihood
  principle.
\newblock \emph{BIT}, 19\penalty0 (2):\penalty0 214--222, 1979{\natexlab{a}}.
\newblock URL \url{https://doi.org/10.1007/BF01930851}.

\bibitem[Larkin(1979{\natexlab{b}})]{Larkin1979a}
F.~M. Larkin.
\newblock Bayesian estimation of zeros of analytic functions.
\newblock Technical report, Queen's University of Kingston. Department of
  Computing and Information Science., 1979{\natexlab{b}}.

\bibitem[Larkin(1979{\natexlab{c}})]{Larkin1979b}
F.~M. Larkin.
\newblock Probabilistic estimation of poles or zeros of functions.
\newblock \emph{J. Approx. Theory}, 27\penalty0 (4):\penalty0 355--371,
  1979{\natexlab{c}}.
\newblock URL \url{https://doi.org/10.1016/0021-9045(79)90124-2}.

\bibitem[Larkin et~al.(1967)Larkin, Brown, Morton, and Bond]{LarkinVideo}
F.~M. Larkin, C.~E. Brown, K.~W. Morton, and P.~Bond.
\newblock Worth a thousand words, 1967.
\newblock URL
  \url{http://www.amara.org/en/videos/7De21CeNlz8b/info/worth-a-thousand-words-1967/}.

\bibitem[Lauritzen(1996)]{Lauritzen1996}
S.~L. Lauritzen.
\newblock \emph{Graphical {Models}}, volume~17 of \emph{Oxford Statistical
  Science Series}.
\newblock The Clarendon Press, Oxford University Press, New York, 1996.
\newblock Oxford Science Publications.

\bibitem[Le~Ma\^{i}tre and Knio(2010)]{LeMaitreKnio2010}
O.~P. Le~Ma\^{i}tre and O.~M. Knio.
\newblock \emph{Spectral Methods for Uncertainty Quantification}.
\newblock Scientific Computation. Springer, New York, 2010.
\newblock URL \url{https://doi.org/10.1007/978-90-481-3520-2}.

\bibitem[Lie et~al.(2018)Lie, Sullivan, and
  Teckentrup]{LieSullivanTeckentrup2018}
H.~C. Lie, T.~J. Sullivan, and A.~L. Teckentrup.
\newblock Random forward models and log-likelihoods in {Bayesian} inverse
  problems.
\newblock \emph{SIAM/ASA J. Uncertain. Quantif.}, 6\penalty0 (4):\penalty0
  1600--1629, 2018.
\newblock URL \url{https://doi.org/10.1137/18M1166523}.

\bibitem[Lie et~al.(2019)Lie, Stuart, and Sullivan]{LieStuartSullivan2019}
H.~C. Lie, A.~M. Stuart, and T.~J. Sullivan.
\newblock Strong convergence rates of probabilistic integrators for ordinary
  differential equations.
\newblock \emph{Stat. Comp.}, 2019.
\newblock To appear. \arXiv{1703.03680}.

\bibitem[Lindley(1956)]{Lindley1956}
D.~V. Lindley.
\newblock On a measure of the information provided by an experiment.
\newblock \emph{Ann. Math. Stat.}, 27:\penalty0 986--1005, 1956.
\newblock URL \url{https://doi.org/10.1214/aoms/1177728069}.

\bibitem[Mahsereci and Hennig(2015)]{Mahsereci2015}
M.~Mahsereci and P.~Hennig.
\newblock Probabilistic line searches for stochastic optimization.
\newblock In \emph{Advances in Neural Information Processing Systems 28}, pages
  181--189, 2015.
\newblock URL
  \url{https://papers.nips.cc/paper/5753-probabilistic-line-searches-for-stochastic-optimization}.

\bibitem[Manber and Norvig(2012)]{Manber2012}
U.~Manber and P.~Norvig.
\newblock The power of the {Apollo} missions in a single {Google} search, 2012.
\newblock URL
  \url{https://search.googleblog.com/2012/08/the-power-of-apollo-missions-in-single.html}.

\bibitem[Marques et~al.(2013)Marques, Bouville, Ribardiere, Santos, and
  Bouatouch]{Marques2013}
R.~Marques, C.~Bouville, M.~Ribardiere, L.~P. Santos, and K.~Bouatouch.
\newblock A spherical {Gaussian} framework for {Bayesian} {Monte} {Carlo}
  rendering of glossy surfaces.
\newblock \emph{IEEE Trans. Vis. and Comp. Graph.}, 19\penalty0 (10):\penalty0
  1619--1632, 2013.
\newblock URL \url{https://doi.org/10.1109/TVCG.2013.79}.

\bibitem[Minka(2000)]{Minka:2000}
T.~Minka.
\newblock Deriving quadrature rules from {Gaussian} processes, 2000.
\newblock URL
  \url{https://www.microsoft.com/en-us/research/publication/deriving-quadrature-rules-gaussian-processes/}.

\bibitem[Mo\v{c}kus(1975)]{Mockus1974}
J.~Mo\v{c}kus.
\newblock On {Bayesian} methods for seeking the extremum.
\newblock In \emph{Optimization Techniques IFIP Technical Conference
  Novosibirsk, July 1--7, 1974. Optimization Techniques 1974}, volume~27 of
  \emph{Lecture Notes in Computer Science}, pages 400--404. Springer, Berlin,
  Heidelberg, 1975.
\newblock URL \url{https://doi.org/10.1007/3-540-07165-2\_55}.

\bibitem[Mo\v{c}kus(1977)]{Mockus1977}
J.~Mo\v{c}kus.
\newblock On {Bayesian} methods for seeking the extremum and their application.
\newblock In \emph{Information {Processing} 77 (Proc.\ IFIP Congr., Toronto,
  Ont., 1977)}, pages 195--200. IFIP Congr. Ser., Vol. 7. North-Holland,
  Amsterdam, 1977.

\bibitem[Mo\v{c}kus(1989)]{Mockus1989}
J.~Mo\v{c}kus.
\newblock \emph{Bayesian approach to global optimization}, volume~37 of
  \emph{Mathematics and its Applications (Soviet Series)}.
\newblock Kluwer Academic Publishers Group, Dordrecht, 1989.
\newblock URL \url{https://doi.org/10.1007/978-94-009-0909-0}.

\bibitem[Murphy(2002)]{Murphy2002}
K.~P. Murphy.
\newblock \emph{Dynamic {Bayesian} networks: representation, inference and
  learning}.
\newblock PhD thesis, University of California, Berkeley, 2002.

\bibitem[Norden et~al.(1978)Norden, Zabotin, {\`E}skin, Grigor{\cprime}ev, and
  Begovatov]{Norden1978}
A.~P. Norden, Y.~I. Zabotin, L.~D. {\`E}skin, S.~V. Grigor{\cprime}ev, and
  E.~A. Begovatov.
\newblock Al{\cprime}bert {Valentinovich} {Sul{\cprime}din} (on the occasion of
  his fiftieth birthday).
\newblock \emph{Izv. Vys\v{s}. U\v{c}ebn. Zaved. Mat.}, 12:\penalty0 3--5,
  1978.

\bibitem[Novak(1988)]{Novak1988}
E.~Novak.
\newblock \emph{Deterministic and {Stochastic} {Error} {Bounds} in {Numerical}
  {Analysis}}, volume 1349 of \emph{Lecture Notes in Mathematics}.
\newblock Springer-Verlag, Berlin, 1988.
\newblock URL \url{https://doi.org/10.1007/BFb0079792}.

\bibitem[Oates et~al.(2017)Oates, Niederer, Lee, Briol, and
  Girolami]{Oates2016}
C.~Oates, S.~Niederer, A.~Lee, F.-X. Briol, and M.~Girolami.
\newblock Probabilistic models for integration error in the assessment of
  functional cardiac models.
\newblock In \emph{Advances in Neural Information Processing Systems 30}, pages
  110--118, 2017.
\newblock URL
  \url{http://papers.nips.cc/paper/6616-probabilistic-models-for-integration-error-in-the-assessment-of-functional-cardiac-models}.

\bibitem[Oates et~al.(2019{\natexlab{a}})Oates, Cockayne, Aykroyd, and
  Girolami]{Oates2017}
C.~J. Oates, J.~Cockayne, R.~G. Aykroyd, and M.~Girolami.
\newblock Bayesian probabilistic numerical methods in time-dependent state
  estimation for industrial hydrocyclone equipment.
\newblock \emph{J. Amer. Stat. Assoc.}, 2019{\natexlab{a}}.
\newblock URL \url{https://doi.org/10.1080/01621459.2019.1574583}.

\bibitem[Oates et~al.(2019{\natexlab{b}})Oates, Cockayne, Prangle, Sullivan,
  and Girolami]{Oates_RICAM}
C.~J. Oates, J.~Cockayne, D.~Prangle, T.~J. Sullivan, and M.~Girolami.
\newblock Optimality criteria for probabilistic numerical methods.
\newblock In \emph{Multivariate Algorithms and Information-Based Complexity,
  Linz, 2018}, 2019{\natexlab{b}}.
\newblock \arXiv{1901.04326}.

\bibitem[O'Hagan(1978)]{OHagan1978}
A.~O'Hagan.
\newblock Curve fitting and optimal design for prediction.
\newblock \emph{J. Roy. Stat. Soc. Ser. B}, 40\penalty0 (1):\penalty0 1--42,
  1978.
\newblock URL \url{https://doi.org/10.1111/j.2517-6161.1978.tb01643.x}.

\bibitem[O'Hagan(1991)]{OHagan1991}
A.~O'Hagan.
\newblock {Bayes}--{Hermite} quadrature.
\newblock \emph{J. Stat. Plann. Inference}, 29\penalty0 (3):\penalty0 245--260,
  1991.
\newblock URL \url{https://doi.org/10.1016/0378-3758(91)90002-V}.

\bibitem[Osborne(2014{\natexlab{a}})]{Osborne2014a}
M.~Osborne.
\newblock {T}\"{u}bingen manifesto: Uncertainty, 2014{\natexlab{a}}.
\newblock URL
  \url{http://probabilistic-numerics.org/2014/08/27/Roundtable-Uncertainty/}.

\bibitem[Osborne(2014{\natexlab{b}})]{Osborne2014b}
M.~Osborne.
\newblock {T}\"{u}bingen manifesto: Probabilistic numerics and probabilistic
  programming, 2014{\natexlab{b}}.
\newblock URL
  \url{http://probabilistic-numerics.org/2014/09/01/Roundtable-ProbNum-ProbProg/}.

\bibitem[Osborne(2014{\natexlab{c}})]{Osborne2014c}
M.~Osborne.
\newblock {T}\"{u}bingen manifesto: Priors and prior work, 2014{\natexlab{c}}.
\newblock URL
  \url{http://probabilistic-numerics.org/2014/08/27/Roundtable-Uncertainty/}.

\bibitem[Osborne(2014{\natexlab{d}})]{Osborne2014d}
M.~Osborne.
\newblock {T}\"{u}bingen manifesto: Community, 2014{\natexlab{d}}.
\newblock URL
  \url{http://probabilistic-numerics.org/2014/09/05/Roundtable-Community/}.

\bibitem[Osborne et~al.(2012{\natexlab{a}})Osborne, Garnett, Ghahramani,
  Duvenaud, Roberts, and Rasmussen]{Osborne2012a}
M.~Osborne, R.~Garnett, Z.~Ghahramani, D.~K. Duvenaud, S.~J. Roberts, and C.~E.
  Rasmussen.
\newblock Active learning of model evidence using {Bayesian} quadrature.
\newblock In \emph{Advances in Neural Information Processing Systems 25}, pages
  46--54, 2012{\natexlab{a}}.
\newblock URL
  \url{https://papers.nips.cc/paper/4657-active-learning-of-model-evidence-using-bayesian-quadrature}.

\bibitem[Osborne et~al.(2012{\natexlab{b}})Osborne, Garnett, Roberts, Hart,
  Aigrain, Gibson, and Aigrain]{Osborne2012}
M.~A. Osborne, R.~Garnett, S.~J. Roberts, C.~Hart, S.~Aigrain, N.~Gibson, and
  S.~Aigrain.
\newblock Bayesian quadrature for ratios.
\newblock In \emph{Proceedings of Artificial Intelligence and Statistics
  (AISTATS)}, 2012{\natexlab{b}}.

\bibitem[Owen(2019)]{Owen2019}
A.~Owen.
\newblock Unreasonable effectiveness of {Monte} {Carlo}.
\newblock \emph{Stat. Sci.}, 2019.

\bibitem[Owhadi(2015)]{Owhadi2015}
H.~Owhadi.
\newblock Bayesian numerical homogenization.
\newblock \emph{Multiscale Model. Simul.}, 13\penalty0 (3):\penalty0 812--828,
  2015.
\newblock URL \url{https://doi.org/10.1137/140974596}.

\bibitem[Owhadi(2017)]{Owhadi2017sirev}
H.~Owhadi.
\newblock Multigrid with rough coefficients and multiresolution operator
  decomposition from hierarchical information games.
\newblock \emph{SIAM Rev.}, 59\penalty0 (1):\penalty0 99--149, 2017.
\newblock URL \url{https://doi.org/10.1137/15M1013894}.

\bibitem[Owhadi and Scovel(2017{\natexlab{a}})]{Owhadi2017general}
H.~Owhadi and C.~Scovel.
\newblock Universal scalable robust solvers from computational information
  games and fast eigenspace adapted multiresolution analysis,
  2017{\natexlab{a}}.
\newblock \arXiv{1703.10761}.

\bibitem[Owhadi and Scovel(2017{\natexlab{b}})]{OwhadiScovel2017}
H.~Owhadi and C.~Scovel.
\newblock Toward {Machine} {Wald}.
\newblock In \emph{Handbook of Uncertainty Quantification}, pages 157--191.
  Springer International Publishing, 2017{\natexlab{b}}.
\newblock URL \url{https://doi.org/10.1007/978-3-319-12385-1\_3}.

\bibitem[Owhadi and Zhang(2017)]{Owhadi2017jcp}
H.~Owhadi and L.~Zhang.
\newblock Gamblets for opening the complexity-bottleneck of implicit schemes
  for hyperbolic and parabolic {ODE}s/{PDE}s with rough coefficients.
\newblock \emph{J. Comp. Phys.}, 347:\penalty0 99--128, 2017.
\newblock URL \url{https://doi.org/10.1016/j.jcp.2017.06.037}.

\bibitem[Owhadi et~al.(2015)Owhadi, Scovel, and
  Sullivan]{Owhadi2015_brittleness}
H.~Owhadi, C.~Scovel, and T.~J. Sullivan.
\newblock Brittleness of {Bayesian} inference under finite information in a
  continuous world.
\newblock \emph{Electron. J. Stat.}, 9\penalty0 (1):\penalty0 1--79, 2015.
\newblock URL \url{https://doi.org/10.1214/15-EJS989}.

\bibitem[Parzen(1970)]{Parzen1970}
E.~Parzen.
\newblock Statistical inference on time series by {RKHS} methods.
\newblock Technical report, Stanford University of California, Department of
  Statistics, 1970.

\bibitem[Paul et~al.(2018)Paul, Chatzilygeroudis, Ciosek, Mouret, Osborne, and
  Whiteson]{Paul2018}
S.~Paul, K.~Chatzilygeroudis, K.~Ciosek, J.-B. Mouret, M.~A. Osborne, and
  S.~Whiteson.
\newblock Alternating optimisation and quadrature for robust control.
\newblock In \emph{The Thirty-Second AAAI Conference on Artificial Intelligence
  (AAAI-18)}, 2018.

\bibitem[Piiroinen(2005)]{Piiroinen2005}
P.~Piiroinen.
\newblock \emph{Statistical {Measurements}, {Experiments} and {Applications}}.
\newblock PhD thesis, University of Helsinki, 2005.

\bibitem[Poincar{\'e}(1896)]{Poincare1896}
H.~Poincar{\'e}.
\newblock \emph{Calcul des {Probabilit{\'e}s}}.
\newblock Georges Carr{\'e}, 1896.

\bibitem[Poincar{\'e}(1912)]{Poincare1912}
H.~Poincar{\'e}.
\newblock \emph{Calcul des {Probabilit{\'e}s}}.
\newblock Gauthier-Villars, second edition, 1912.

\bibitem[Pr\"{u}her et~al.(2018)Pr\"{u}her, Karvonen, Oates, Straka, and
  S\"{a}rkk\"{a}]{Pruher2018}
J.~Pr\"{u}her, T.~Karvonen, C.~J. Oates, O.~Straka, and S.~S\"{a}rkk\"{a}.
\newblock Improved calibration of numerical integration error in sigma-point
  filters, 2018.
\newblock \arXiv{1811.11474}.

\bibitem[{Queen's University at Kingston}(11 Feb.\ 1982)]{Queens1982}
{Queen's University at Kingston}.
\newblock Frederick {Michael} {Larkin} (1936--1982), 11 Feb.\ 1982.
\newblock URL
  \url{https://grahamlarkin.files.wordpress.com/2018/12/fmlarkin\_obit.pdf}.

\bibitem[Rasmussen and Ghahramani(2003)]{RasmussenGhahramani:2003}
C.~E. Rasmussen and Z.~Ghahramani.
\newblock Bayesian {Monte} {Carlo}.
\newblock In \emph{Advances in Neural Information Processing Systems 16}, pages
  505--512, 2003.
\newblock URL \url{http://papers.nips.cc/paper/2150-bayesian-monte-carlo}.

\bibitem[Ritter(2000)]{Ritter2000}
K.~Ritter.
\newblock \emph{Average-{C}ase {A}nalysis of {N}umerical {P}roblems}, volume
  1733 of \emph{Lecture Notes in Mathematics}.
\newblock Springer-Verlag, Berlin, 2000.
\newblock URL \url{https://doi.org/10.1007/BFb0103934}.

\bibitem[Robins and van~der Vaart(2006)]{Robins2006}
J.~Robins and A.~van~der Vaart.
\newblock Adaptive nonparametric confidence sets.
\newblock \emph{Ann. Stat.}, 34\penalty0 (1):\penalty0 229--253, 2006.
\newblock URL \url{https://doi.org/10.1214/009053605000000877}.

\bibitem[Sacks and Ylvisaker(1968)]{SacksYlvisaker1968}
J.~Sacks and D.~Ylvisaker.
\newblock Designs for regression problems with correlated errors; many
  parameters.
\newblock \emph{Ann. Math. Stat.}, 39:\penalty0 49--69, 1968.
\newblock URL \url{https://doi.org/10.1214/aoms/1177698504}.

\bibitem[Sacks and Ylvisaker(1970{\natexlab{a}})]{SacksYlvisaker1970a}
J.~Sacks and D.~Ylvisaker.
\newblock Designs for regression problems with correlated errors. {III}.
\newblock \emph{Ann. Math. Stat.}, 41:\penalty0 2057--2074, 1970{\natexlab{a}}.
\newblock URL \url{https://doi.org/10.1214/aoms/1177696705}.

\bibitem[Sacks and Ylvisaker(1970{\natexlab{b}})]{SacksYlvisaker1970b}
J.~Sacks and D.~Ylvisaker.
\newblock Statistical designs and integral approximation.
\newblock In \emph{Proc. {Twelfth} {Biennial} {Sem.} {Canad.} {Math.} {C}ongr.
  on {T}ime {S}eries and {S}tochastic {P}rocesses; {C}onvexity and
  {C}ombinatorics ({V}ancouver, {B}.{C}., 1969)}, pages 115--136. Canad. Math.
  Congr., Montreal, Que., 1970{\natexlab{b}}.

\bibitem[Sacks and Ylvisaker(1966)]{SacksYlvisaker1966}
J.~Sacks and N.~D. Ylvisaker.
\newblock Designs for regression problems with correlated errors.
\newblock \emph{Ann. Math. Stat.}, 37:\penalty0 66--89, 1966.
\newblock URL \url{https://doi.org/10.1214/aoms/1177699599}.

\bibitem[Sacks et~al.(1989)Sacks, Welch, Mitchell, and
  Wynn]{SacksWelchMitchellWynn1989}
J.~Sacks, W.~J. Welch, T.~J. Mitchell, and H.~P. Wynn.
\newblock Design and analysis of computer experiments.
\newblock \emph{Stat. Sci.}, 4\penalty0 (4):\penalty0 409--435, 1989.
\newblock URL \url{https://doi.org/10.1214/ss/1177012413}.

\bibitem[Sard(1949)]{Sard1949}
A.~Sard.
\newblock Best approximate integration formulas; best approximation formulas.
\newblock \emph{Amer. J. Math.}, 71:\penalty0 80--91, 1949.
\newblock URL \url{https://doi.org/10.2307/2372095}.

\bibitem[Sard(1963)]{Sard1963}
A.~Sard.
\newblock \emph{Linear {Approximation}}.
\newblock Number~9 in Mathematical Surveys. American Mathematical Society,
  Providence, RI, 1963.
\newblock URL \url{https://doi.org/10.1090/surv/009}.

\bibitem[S{\"a}rkk{\"a} et~al.(2016)S{\"a}rkk{\"a}, Hartikainen, Svensson, and
  Sandblom]{Saerkkae2015}
S.~S{\"a}rkk{\"a}, J.~Hartikainen, L.~Svensson, and F.~Sandblom.
\newblock On the relation between {Gaussian} process quadratures and
  sigma-point methods.
\newblock \emph{J. Adv. Inf. Fusion}, 11\penalty0 (1):\penalty0 31--46, 2016.

\bibitem[Sch{\"a}fer et~al.(2017)Sch{\"a}fer, Sullivan, and
  Owhadi]{Schaefer2017}
F.~Sch{\"a}fer, T.~J. Sullivan, and H.~Owhadi.
\newblock Compression, inversion, and approximate {PCA} of dense kernel
  matrices at near-linear computational complexity, 2017.
\newblock \arXiv{1706.02205}.

\bibitem[Schober et~al.(2014)Schober, Duvenaud, and Hennig]{Schober:2014wt}
M.~Schober, D.~K. Duvenaud, and P.~Hennig.
\newblock Probabilistic {ODE} solvers with {Runge}--{Kutta} means.
\newblock In \emph{Advances in Neural Information Processing Systems 27}, 2014.
\newblock URL
  \url{https://papers.nips.cc/paper/5451-probabilistic-ode-solvers-with-runge-kutta-means}.

\bibitem[Schober et~al.(2018)Schober, S{\"a}rkk{\"a}, and
  Hennig]{Schober:2016uh}
M.~Schober, S.~S{\"a}rkk{\"a}, and P.~Hennig.
\newblock A probabilistic model for the numerical solution of initial value
  problems.
\newblock \emph{Stat. Comp.}, 29\penalty0 (1):\penalty0 99--122, 2018.
\newblock URL \url{https://doi.org/10.1007/s11222-017-9798-7}.

\bibitem[Schoenberg(1965)]{Schoenberg1965}
I.~J. Schoenberg.
\newblock On monosplines of least deviation and best quadrature formulae.
\newblock \emph{J. Soc. Indust. Appl. Math. Ser. B Numer. Anal.}, 2\penalty0
  (1):\penalty0 144--170, 1965.
\newblock URL \url{https://doi.org/10.1137/0702012}.

\bibitem[Schoenberg(1966)]{Schoenberg1966}
I.~J. Schoenberg.
\newblock On monosplines of least square deviation and best quadrature
  formulae. {II}.
\newblock \emph{SIAM J. Numer. Anal.}, 3\penalty0 (2):\penalty0 321--328, 1966.
\newblock URL \url{https://doi.org/10.1137/0703025}.

\bibitem[\'{S}cibior et~al.(2015)\'{S}cibior, Ghahramani, and
  Gordon]{Scibior2015}
A.~\'{S}cibior, Z.~Ghahramani, and A.~Gordon.
\newblock Practical probabilistic programming with monads.
\newblock \emph{ACM SIGPLAN Notices}, 50\penalty0 (12):\penalty0 165--176,
  2015.
\newblock URL \url{https://doi.org/10.1145/2804302.2804317}.

\bibitem[Sirignano and Spiliopoulos(2018)]{Sirignano2018}
J.~Sirignano and K.~Spiliopoulos.
\newblock {DGM}: {A} deep learning algorithm for solving partial differential
  equations.
\newblock \emph{J. Comp. Phys.}, 375:\penalty0 1339--1364, 2018.
\newblock URL \url{https://doi.org/10.1016/j.jcp.2018.08.029}.

\bibitem[Skilling(1992)]{Skilling1992}
J.~Skilling.
\newblock Bayesian solution of ordinary differential equations.
\newblock In \emph{Maximum Entropy and Bayesian Methods}, pages 23--37.
  Springer, 1992.
\newblock URL \url{https://doi.org/10.1007/978-94-017-2219-3}.

\bibitem[Smale(1985)]{Smale1985}
S.~Smale.
\newblock On the efficiency of algorithms of analysis.
\newblock \emph{Bull. Amer. Math. Soc. (N.S.)}, 13\penalty0 (2):\penalty0
  87--121, 1985.
\newblock URL \url{https://doi.org/10.1090/S0273-0979-1985-15391-1}.

\bibitem[Smith(2014)]{Smith:2014}
R.~C. Smith.
\newblock \emph{Uncertainty {Quantification}: {Theory}, {Implementation}, and
  {Applications}}, volume~12 of \emph{Computational Science \& Engineering}.
\newblock Society for Industrial and Applied Mathematics (SIAM), Philadelphia,
  PA, 2014.

\bibitem[Snoek et~al.(2012)Snoek, Larochelle, and Adams]{Snoek2012}
J.~Snoek, H.~Larochelle, and R.~P. Adams.
\newblock Practical {B}ayesian optimization of machine learning algorithms.
\newblock In \emph{Advances in Neural Information Processing Systems}, pages
  2951--2959, 2012.
\newblock URL
  \url{https://papers.nips.cc/paper/4522-practical-bayesian-optimization-of-machine-learning-algorithms}.

\bibitem[Stuart(2010)]{Stuart:2010}
A.~M. Stuart.
\newblock Inverse problems: {A} {Bayesian} perspective.
\newblock \emph{Acta Numer.}, 19:\penalty0 451--559, 2010.
\newblock URL \url{https://doi.org/10.1017/S0962492910000061}.

\bibitem[Stuart and Teckentrup(2018)]{StuartTeckentrup2018}
A.~M. Stuart and A.~L. Teckentrup.
\newblock Posterior consistency for {Gaussian} process approximations of
  {Bayesian} posterior distributions.
\newblock \emph{Math. Comp.}, 87\penalty0 (310):\penalty0 721--753, 2018.
\newblock URL \url{https://doi.org/10.1090/mcom/3244}.

\bibitem[Sudakov(1959)]{Sudakov1959}
V.~N. Sudakov.
\newblock Linear sets with quasi-invariant measure.
\newblock \emph{Dokl. Akad. Nauk SSSR}, 127:\penalty0 524--525, 1959.

\bibitem[Sul{\cprime}din(1959)]{Suldin1959}
A.~V. Sul{\cprime}din.
\newblock Wiener measure and its applications to approximation methods. {I}.
\newblock \emph{Izv. Vys\v{s}. U\v{c}ebn. Zaved. Mat.}, 6\penalty0
  (13):\penalty0 145--158, 1959.

\bibitem[Sul{\cprime}din(1960)]{Suldin1960}
A.~V. Sul{\cprime}din.
\newblock Wiener measure and its applications to approximation methods. {II}.
\newblock \emph{Izv. Vys\v{s}. U\v{c}ebn. Zaved. Mat.}, 5\penalty0
  (18):\penalty0 165--179, 1960.

\bibitem[Sul{\cprime}din(1963{\natexlab{a}})]{Suldin1963a}
A.~V. Sul{\cprime}din.
\newblock The method of regression in the theory of approximation.
\newblock \emph{Kazan. Gos. Univ. U\v{c}en. Zap.}, 123\penalty0 (kn.
  6):\penalty0 3--35, 1963{\natexlab{a}}.

\bibitem[Sul{\cprime}din(1963{\natexlab{b}})]{Suldin1963b}
A.~V. Sul{\cprime}din.
\newblock On the distribution of the functional {$\int_{0}^{1} x^{2}(t) \, \rd
  t$} where {$x(t)$} represents a certain {Gaussian} process.
\newblock In \emph{Kazan {State} {Univ.} {Sci.} {Survey} {Conf.} 1962
  ({Russian})}, pages 80--82. Izdat. Kazan. Univ., Kazan, 1963{\natexlab{b}}.

\bibitem[Sul{\cprime}din(1963{\natexlab{c}})]{Suldin1963c}
A.~V. Sul{\cprime}din.
\newblock The solution of equations by the method of conditional mean values.
\newblock In \emph{Kazan {S}tate {U}niv. {S}ci. {S}urvey {C}onf. 1962
  ({R}ussian)}, pages 85--87. Izdat. Kazan. Univ., Kazan, 1963{\natexlab{c}}.

\bibitem[Sul{\cprime}din(1968)]{Suldin1968}
A.~V. Sul{\cprime}din.
\newblock Curves and operators in a {Hilbert} space.
\newblock \emph{Kazan. Gos. Univ. U\v{c}en. Zap.}, 128\penalty0 (2):\penalty0
  15--47, 1968.

\bibitem[Sul{\cprime}din et~al.(1969)Sul{\cprime}din, Zabotin, and
  Semenihina]{Suldin1969}
A.~V. Sul{\cprime}din, V.~I. Zabotin, and N.~P. Semenihina.
\newblock Certain operators in {Hilbert} space.
\newblock \emph{Kazan. Gos. Univ. U\v{c}en. Zap.}, 129\penalty0 (4):\penalty0
  90--95, 1969.

\bibitem[Sullivan(2015)]{Sullivan:2015}
T.~J. Sullivan.
\newblock \emph{Introduction to {Uncertainty} {Quantification}}, volume~63 of
  \emph{Texts in Applied Mathematics}.
\newblock Springer, 2015.
\newblock URL \url{https://doi.org/10.1007/978-3-319-23395-6}.

\bibitem[Teymur et~al.(2016)Teymur, Zygalakis, and Calderhead]{Teymur2016}
O.~Teymur, K.~Zygalakis, and B.~Calderhead.
\newblock Probabilistic linear multistep methods.
\newblock In \emph{Advances in Neural Information Processing Systems 29}, 2016.
\newblock URL
  \url{https://papers.nips.cc/paper/6356-probabilistic-linear-multistep-methods}.

\bibitem[Teymur et~al.(2018)Teymur, Lie, Sullivan, and Calderhead]{Teymur2018}
O.~Teymur, H.~C. Lie, T.~J. Sullivan, and B.~Calderhead.
\newblock Implicit probabilistic integrators for {ODE}s.
\newblock In \emph{32nd Conference on Neural Information Processing Systems
  (NeurIPS 2018)}, 2018.
\newblock URL
  \url{http://papers.nips.cc/paper/7955-implicit-probabilistic-integrators-for-odes}.

\bibitem[{The MathWorks Inc.}()]{BayesOptMatlab}
{The MathWorks Inc.}
\newblock Bayesian optimization algorithm.
\newblock URL
  \url{https://uk.mathworks.com/help/stats/bayesian-optimization-algorithm.html}.
\newblock Accessed December 2018.

\bibitem[Tienari(1970)]{Tienari1970}
M.~Tienari.
\newblock A statistical model of roundoff error for varying length
  floating-point arithmetic.
\newblock \emph{Nordisk Tidskr. Informationsbehandling (BIT)}, 10:\penalty0
  355--365, 1970.
\newblock URL \url{https://doi.org/10.1007/BF01934204}.

\bibitem[Traub and Wo\'{z}niakowsi(1980)]{Traub1980}
J.~F. Traub and H.~Wo\'{z}niakowsi.
\newblock \emph{A {General} {Theory} of {Optimal} {Algorithms}}.
\newblock ACM Monograph Series. Academic Press, Inc. [Harcourt Brace
  Jovanovich, Publishers], New York-London, 1980.

\bibitem[Traub et~al.(1983)Traub, Wasilkowski, and Wo{\'z}niakowski]{Traub1983}
J.~F. Traub, G.~W. Wasilkowski, and H.~Wo{\'z}niakowski.
\newblock \emph{Information, {U}ncertainty, {C}omplexity}.
\newblock Addison-Wesley Publishing Company, Advanced Book Program, Reading,
  MA, 1983.

\bibitem[Trefethen(2008)]{Trefethen2008}
L.~N. Trefethen.
\newblock Is {Gauss} quadrature better than {Clenshaw}--{Curtis}?
\newblock \emph{SIAM Rev.}, 50\penalty0 (1):\penalty0 67--87, 2008.
\newblock URL \url{http://dx.doi.org/10.1137/060659831}.

\bibitem[Tronarp et~al.(2019)Tronarp, Kersting, S{\"a}rkk{\"a}, and
  Hennig]{Tronarp2019}
F.~Tronarp, H.~Kersting, S.~S{\"a}rkk{\"a}, and P.~Hennig.
\newblock Probabilistic solutions to ordinary differential equations as
  non-linear {Bayesian} filtering: {A} new perspective, 2019.
\newblock \arXiv{1810.03440}.

\bibitem[{U.S. Department of Energy}(2009)]{GrandChallenges2009}
{U.S. Department of Energy}.
\newblock \emph{Scientific {Grand} {Challenges} for {National} {Security}:
  {The} {Role} of {Computing} at the {Extreme} {Scale}}.
\newblock 2009.

\bibitem[von Neumann and Goldstine(1947)]{VonNeumannGoldstine1947}
J.~von Neumann and H.~H. Goldstine.
\newblock Numerical inverting of matrices of high order.
\newblock \emph{Bull. Amer. Math. Soc.}, 53:\penalty0 1021--1099, 1947.
\newblock URL \url{https://doi.org/10.1090/S0002-9904-1947-08909-6}.

\bibitem[Wang et~al.(2018)Wang, Cockayne, and Oates]{Wang:2018}
J.~Wang, J.~Cockayne, and C.~Oates.
\newblock On the {Bayesian} solution of differential equations.
\newblock In \emph{Bayesian Inference and Maximum Entropy Methods in Science
  and Engineering (MaxEnt 2018)}, 2018.

\bibitem[Wo\'{z}niakowski(1985)]{Wozniakowski1985}
H.~Wo\'{z}niakowski.
\newblock A survey of information-based complexity.
\newblock \emph{J. Complexity}, 1\penalty0 (1):\penalty0 11--44, 1985.
\newblock URL \url{https://doi.org/10.1016/0885-064X(85)90020-2}.

\bibitem[Xi et~al.(2018)Xi, Briol, and Girolami]{Xi2018}
X.~Xi, F.-X. Briol, and M.~Girolami.
\newblock Bayesian quadrature for multiple related integrals.
\newblock In \emph{Proceedings of the 35th International Conference on Machine
  Learning}, volume~80, pages 5373--5382, 2018.
\newblock URL \url{http://proceedings.mlr.press/v80/xi18a/xi18a.pdf}.

\bibitem[Yoo and Owhadi(2019)]{Yoo2019}
G.~R. Yoo and H.~Owhadi.
\newblock De-noising by thresholding operator adapted wavelets.
\newblock \emph{Stat. Comp.}, 2019.
\newblock To appear. \arXiv{1805.10736}.

\bibitem[Zabotin et~al.(1996)Zabotin, Zamov, Aksent{\cprime}ev, and
  Zemtseva]{Zabotin1996}
Y.~I. Zabotin, N.~K. Zamov, L.~A. Aksent{\cprime}ev, and T.~N. Zemtseva.
\newblock Al{\cprime}bert {Valentinovich} {Sul{\cprime}din} (obituary).
\newblock \emph{Izv. Vys\v{s}. U\v{c}ebn. Zaved. Mat.}, 2\penalty0 (84), 1996.

\bibitem[Zellner(1988)]{Zellner:1988}
A.~Zellner.
\newblock Optimal information processing and {Bayes}'s theorem.
\newblock \emph{Amer. Stat.}, 42\penalty0 (4):\penalty0 278--284, 1988.
\newblock URL \url{https://doi.org/10.2307/2685143}.

\end{thebibliography}

\end{document}